\documentclass[a4wide,twoside,12pt]{article}

\usepackage{a4wide}
\usepackage{graphicx}
\usepackage{amstext, amsmath, amssymb, amsfonts, amsbsy}

\usepackage{latexsym}
\usepackage{booktabs}
\usepackage{url}
\usepackage{import}
\usepackage{fancyvrb}
\usepackage{stmaryrd}
\usepackage{pdfsync}
\usepackage{algorithm}
\usepackage{tikz}
\usepackage{pgfplots}
\usepackage{pgfplotstable}
\usepackage{subfig}
\usepackage{tabularx}
\usepackage{snapshot}
\usepackage{algpseudocode}

\usepackage[shadow]{todonotes}
\usepackage[numbers]{natbib}

\newtheorem{remark}{Remark}[section] 

\numberwithin{equation}{section}
\numberwithin{figure}{section}
\numberwithin{table}{section}

\newcommand{\R}{\mathbb{R}}
\newcommand{\N}{\mathbb{N}}

\newcommand{\foralls}{\forall\,}

\newcommand{\mesh}{\mathcal{T}_h}

\DefineVerbatimEnvironment{code}{Verbatim}{frame=single,rulecolor=\color{blue}}

\newcommand{\mcF}{\mathcal{F}}

\newcommand{\restr}[2]{ \left. #1 \right|_{#2}}

\newcommand{\jump}[1]{\left\llbracket #1\right\rrbracket}
\newcommand{\norm}[1]{\left\| #1\right\|}

\title{\bf A CutFEM method for Stefan-Signorini problems with application in pulsed laser ablation. \\
}
\author{Susanne Claus${\footnote{ClausS2@cardiff.ac.uk}}$ , Samuel Bigot, Pierre Kerfriden${\footnote{pierre.kerfriden@gmail.com}}$\\
\\
Cardiff University, School of Engineering, \\ The Parade, CF243AA Cardiff, United Kingdom. }

\begin{document}
\maketitle

\begin{abstract}
  \noindent 
In this article, we develop a cut finite element method for one-phase Stefan problems with applications in laser manufacturing. The geometry of the workpiece is represented implicitly via a level set function. Material above the melting/vaporisation temperature is represented by a fictitious gas phase. The moving interface between the workpiece and the fictitious gas phase may cut arbitrarily through the elements of the finite element mesh, which remains fixed throughout the simulation, thereby circumventing the need for cumbersome re-meshing operations. The primal/dual formulation of the linear one-phase Stefan problem is recast into a primal non-linear formulation using a Nitsche-type approach, which avoids the difficulty of constructing inf-sup stable primal/dual pairs. Through the careful derivation of stabilisation terms, we show that the proposed Stefan-Signorini-Nitsche CutFEM method remains stable independently of the cut location. In addition, we obtain optimal convergence with respect to space and time refinement. Several 2D and 3D examples are proposed, highlighting the robustness and flexibility of the algorithm, together with its relevance to the field of micro-manufacturing.
\end{abstract}

CutFEM, Stefan problem, Stefan-Signorini-Nitsche formulation, pulsed laser ablation. 

\section{Introduction}

The simulation of phase changes requires tracking the evolution of solid/liquid and liquid/gas interfaces, which is numerically challenging. In the context of the finite element method (FEM), two main approaches for interface tracking can be distinguished. The first family of approaches smooths the transition between phases, allowing for the existence of a mushy region in space where both phases coexist (\textit{i.e.} enthalpy method \cite{Voller1987,Date1992,Belhamadia2012,Danaila2014}, phase field method \cite{Steinbach2009,Zhao2017}). The width of this region may be thought of as a trade-off between computational cost, which is lower for fatter transition zones, and modelling accuracy, whereby the``true" model corresponds to an infinitely thin transition zone. The second approach describes the interface between phases as a sharp surface in 3D or a line in 2D. Although this may seem to be the ``natural" approach to interface tracking, the sharp interface approach is difficult to handle within a finite element context. Indeed, the mesh either needs to conform to this interface, leading to a class of moving mesh algorithms such as ALE, or special finite element methods need to be developed so as to allow the interface to \textit{cut through} the element. The latter family of methods are the so-called implicit boundary methods (see for instance \cite{Moes1999,bordasmoran2006,Bordas2010,Fries2011,HaHa04,BurmanClausHansboEtAl2014}), which are of prime interest in this paper.

The XFEM method was proposed in \cite{Moes1999}, and relies on a partition-of-unity enrichment to represent embedded kinks and discontinuities. The XFEM method has been applied to the simulation of two-phase Stefan problems in \textit{e.g.} \cite{Merle2002,Chessa2002,Zabaras2006,duddu2008combined,Salvatori2009,Bernauer2011,Cosimo2013, Martin2016, Jahn2016a}. In this case, the interface between solid and liquid moves through a regular background mesh, which may be refined around the interface for accuracy purposes, but does not need to conform to it. XFEM methods offer increased efficiency and robustness as they do not require the cumbersome re-meshing operations that are used in mesh-moving algorithm to prevent the development of prohibitively large mesh deformations. As an alternative to XFEM, the CutFEM approach \cite{HansboHansbo2002,HaHa04,BurmanClausHansboEtAl2014} also enriches elements in order to allow for the representation of embedded discontinuities. However, the enrichment is obtained by an overlapping domain-decomposition strategy (\textit{i.e.} a ``fictitious domain" approach). The strength of CutFEM lies in its stability, which is provided by the so-called ``ghost-penalty"  regularisation \cite{Burman2010}. As far as we are aware, CutFEM has never been applied to Stefan problems.

In this paper, we are interested in a particular subclass of phase change problems with sharp interface representation: the one-phase Stefan problem. In this particular setting, only one of the phases is represented and the other phase is replaced by a fictitious material with zero specific heat. Subsequently, the fictitious phase does not contribute to the energy balance. The interface between the represented phase and the fictitious phase is moved so that the flux at the boundary of the represented domain is balanced by a latent heat term. It is possible to include a non-zero energy flux applied locally at the phase change interface. This is routinely done in laser manufacturing to simulate the irradiation of the ablated material (see \cite{Storti1995}). The boundary conditions of the one-phase Stefan problem are ambiguous and can be treated mathematically using the same tools that are used to formulate unilateral contact in solid mechanics. Mathematical considerations relative to the Stefan-Signorini problem can be found in \cite{Friedman1984, Jiang1985, Ton1994}. 

Very few implicit boundary methods have been developed for one-phase Stefan problems. One exception is the elegant Stefan-Signorini formulation proposed in \cite{Narimanyan2009} for the simulation of thermal plasma cutting, associated with an implicit representation of the domain boundary through the evolution of a level-set function. However, unilateral contact problems have been extensively studied in XFEM \cite{dolbowmoes2001,khoeinikbakht2006,ribeaucourtbaiettodubourg2007,elguedjgravouil2007,gravouilpierres2011,muellerhoeppe2012} and CutFEM \cite{Claus2018,Burman2017}. Typical embedded interface formulations of contact laws include the penalty method, the Lagrange multiplier approach, and the Nitsche-contact formulation, which was recently proposed in \cite{choulyhild2013}. Lagrange multiplier approaches are usually solved by either the augmented Lagrangian algorithm \cite{Tur2015, Navarro-Jimenez2017, Burman2016}, the Uzawa algorithm or the LaTIn approach \cite{kerfridenallix2009,allixkerfriden2010,Claus2018}, all of them being some form of proximal algorithms. Alternatively, the Nitsche-contact formulation reformulates the KKT primal/dual contact problem as a purely primal nonlinear problem that can be solved by Newton algorithms \cite{Burman2017a, choulyhild2013, Chouly2015,Chouly2017}. The Nitsche-contact algorithm promises consistency, whilst circumventing the cumbersome choice of an inf-sup-stable pair for the primal and dual finite element spaces.

In this paper, we present the first CutFEM algorithm for phase-change problems with sharp embedded representation of moving interfaces. The proposed algorithm is highly flexible owing to the implicit description of the domain geometry by a level-set function. The most novel part of the algorithm resides in rewriting the primal/dual condition associated with the interface of the one-phase Stefan-Signorini problem using a dedicated Nitsche reformulation, inspired by \cite{Burman2017a, choulyhild2013} and presented in Section 2 and 3. Noteworthily, we provide the expression of the tangent operator required to deploy the Newton algorithm. In addition, our derivation of the Nitsche-Signorini formulation of the Stefan problem departs from those proposed in \cite{Chouly2017,Burman2017a}, which provides new insights into this emerging approach. Consistently with the CutFEM paradigm, cut elements are regularised using ghost penalties \cite{Burman2010}, which are carefully adapted to the context of the one-phase Stefan-Signorini problem. We pay particular attention to the mathematical scaling of all stabilisation terms so as to obtain the optimal trade-off between stability and accuracy. This is described in Section 4 of the paper. Our numerical time integration procedure is relatively classical: we use an implicit Euler algorithm to solve the unsteady temperature equation. At every time step, the interface is moved by advecting the level-set function using the velocity field delivered by the Nitsche-Signorini algorithm and extended to the entire domain using the fast-marching method \cite{Sethian1999}. The advection of the level-set function is stabilised by the SUPG method \cite{Brooks1982}. This is detailed in Section 5.
An alternative approach where the authors use an upwind finite-difference scheme to update the level-set can be found in \cite{duddu2014numerical}. 

Our formal developments are accompanied by a high-performance computer implementation. The core of our implementation is the finite element C++/Python library FEniCS \cite{fenics2015,Hale2017}, which, in particular, proposes a range of high-level tools to rapidly developed finite element solvers. The CutFEM C++ library, LibCutFEM, partially described in \cite{BurmanClausHansboEtAl2014}, defines additional tools that are specific to unfitted finite elements. This library forms the basis for the CutFEM Stefan-Signorini code used in this article.

Several 2D and 3D examples are presented in Section 6, with particular relevance to engineers interested in the simulation of laser micro-milling. We first derive a new 2D manufactured analytical solution for the one-phase Stefan-Signorini problem, which we use to validate our numerical algorithm and show optimal convergence. In particular, we show that the convergence of the proposed algorithm is optimal: order two in space and order one in time. We then move to the 2D simulation of a pulsed laser ablation process, where we compare the effect of different pulse frequencies onto the finishing quality of the ablated surface. Our first 3D example is the simulation of a laser drilling operation, where the laser irradiates the material. Finally, we present a 3D example of laser milling, where material is removed by the laser through a complex path, following a layer-by-layer removal strategy.

\section{The Stefan-Signorini problem}
\label{sec:tfstokes}
\begin{figure}[t!]
  \begin{center}
    \includegraphics[width=0.7\textwidth]{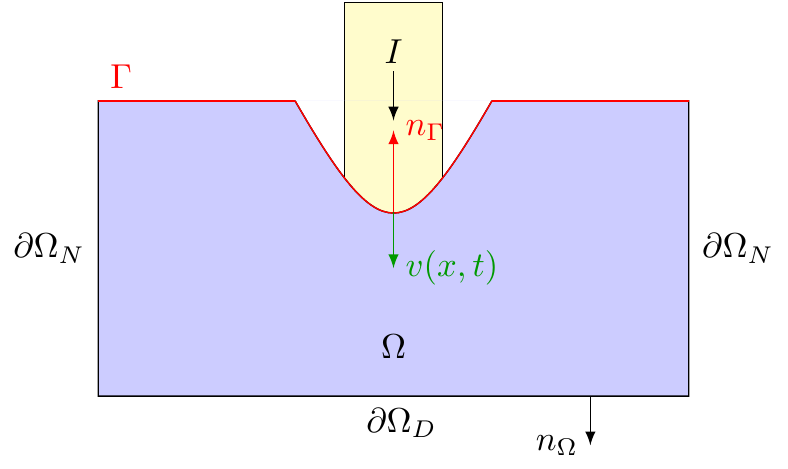}
  \end{center}
  \caption{Schematics of a one-phase Stefan problem with energy flux $I$ heating material $\Omega$ at $\Gamma$.}
\label{fig: laser cutting schematics}
\end{figure}

Let us assume that an energy flux $I(x,t)$ (laser beam) heats a piece of material occupying domain $\Omega(t)$ in $\mathbb{R}^d$ ($d=2$ or $3$) in the time interval $t \in [t_0,t_f]$. Here, $t_0$ is the initial time and $t_f$ is the final time. The boundary of $\Omega(t)$, denoted by $\partial \Omega(t)$, is decomposed into a Dirichlet part $\partial \Omega_D$ and a Neumann part $\partial \Omega_N$ as well as a moving boundary $\Gamma(t)$ (see Figure~\ref{fig: laser cutting schematics}). The boundary $\Gamma(t)$ can be interpreted as an interface between a heated material and the fictitious fluid and gas phase, which will not be represented explicitly. Instead we assume that material above the melting temperature is instantaneously removed, latent heat being consumed in the process. \\
The Stefan-Signorini problem describing the heat conduction and removal of material can then be formulated as:  For all $t \in [t_0,t_f]$, find the temperature  $T: \Omega(t) \rightarrow \mathbb{R}$ such that
  \begin{equation}
  \begin{aligned}
\rho \, c \, \frac{\partial T}{\partial t} - k \, \Delta T &= f \quad &\mbox{ in } \Omega(t) 
\label{equ: parabolic temperature equation}
  \end{aligned}
  \end{equation}
with boundary conditions 
  \begin{equation}
  \begin{aligned}
T &= T_D \quad &\mbox{ on } \partial \Omega_D(t), \\ 
 k \, \nabla T \cdot n_{\Omega} &= q_N   &\mbox{ on } \partial \Omega_N(t) \, ,
\label{equ: strong boundary conditions}
  \end{aligned}
  \end{equation}
together with the evolution equation for interface $\Gamma(t)$ between the fictitious material and the heated material
  \begin{equation} 
I(x,t)  \cdot n_{\Gamma} - k \, \nabla T \cdot n_{\Gamma}  = -\rho \, L \, (v(x,t) \cdot n_{\Gamma}) \quad \mbox{ on } \Gamma(t)
\label{equ: latent heat interface}
  \end{equation}
and the associated Signorini conditions 
  \begin{equation}
  \begin{aligned}
k \, \nabla T \cdot n_{\Gamma} - I \cdot  n_{\Gamma}    &\leq 0 &\mbox{ on } \Gamma(t), \\ 
T - T_m &\leq 0 &\mbox{ on } \Gamma(t), \\
(k \, \nabla T \cdot n_{\Gamma} - I \cdot  n_{\Gamma} ) \perp (T-&T_m)  & \mbox{ on } \Gamma(t),
  \end{aligned}
\label{equ: strong signorini conditions}
  \end{equation}
and the initial condition 
  \begin{equation}
  \begin{aligned}
T(x,t_0) &= T_0 &\mbox{ in } \Omega(t_0),
  \end{aligned}
    \label{eq:initial condition}
  \end{equation}
where $T_0$ is a specified initial temperature with $T_0<T_m$.
Here, $\rho$ is the mass density, $c$ is the heat capacity, $k$ is the thermal conductivity, $f$ is a volumetric heat source, $T_D$ is a given temperature, $q_N$ is a given heat flux, $L$ is the latent heat, $v$ is the velocity of the phase boundary $\Gamma(t)$, $T_m$ is the melting temperature,  and $I$ is the prescribed heat flux given by the profile
\begin{equation}
I(x,t) = A(x,t) \, e_{ray}(x,t),
\end{equation}
where $A$ is the amplitude of the laser beam and  $e_{ray}$ is the direction of the beam.
\begin{remark}
The Signorini conditions \eqref{equ: strong signorini conditions} ensure that material is only removed if it reaches melting temperature and that material is only removed and not added. It can be interpreted as enforcing either \\ 
\textbf{(a)} $T=T_m$ and $ \rho \, L \, v(x,t) \cdot n_{\Gamma} \leq 0$, i.e. the material is heated to melting temperature and material is removed in the normal direction to the interface with speed $v(x,t)$ 
\\
or
\\ 
\textbf{(b)}  $T<T_m$ and $ \rho \, L \, v(x,t) \cdot n_{\Gamma} =0$, i.e. the material does not reach melting temperature and hence no material is removed. 
\end{remark}

\section{Signorini-Nitsche reformulation of the one-phase Stefan problem}
\subsection{Primal/dual weak formulation of the Stefan-Signorini problem}
The weak formulation of the Stefan-Signorini problem reads: For all $t \in [t_0,t_f]$, find $T \in H^1_D(\Omega(t))$ such that, for all $\delta T \in H^1_0(\Omega(t))$, 
\begin{eqnarray}
\rho \, c \, 
\left(\frac{\partial T}{\partial t}, \, \delta T \right)_{\Omega(t)} + \left( k \, \nabla T, \,
\nabla \delta T \right)_{\Omega(t)} =  \left(f , \,
 \delta T \right)_{\Omega(t)} 
+ \left(k \nabla T \cdot n_{\Gamma} , \,
\delta T \right)_{\Gamma(t)} + \left(q_N , \,
 \delta T \right)_{\partial \Omega_N(t)} .
\label{eq:stefan-weak}
\end{eqnarray}
Together with \eqref{equ: latent heat interface}, \eqref{equ: strong signorini conditions} and \eqref{eq:initial condition}. Here, 
\begin{align}
H^1_D(\Omega(t)) &= \{ T \in H^1(\Omega(t)): T = T_D \mbox{ on } \partial \Omega_D(t) \}, \nonumber \\
H^1_0(\Omega(t)) &= \{ T \in H^1(\Omega(t)): T = 0 \mbox{ on } \partial \Omega_D(t) \}.
\end{align}
In order to facilitate the treatment of Signorini conditions \eqref{equ: strong signorini conditions}, we introduce a slack variable 
\begin{equation}
\sigma : = k \, \nabla T \cdot n_{\Gamma} - I \cdot  n_{\Gamma},
\label{equ: slack variable}
\end{equation}
which yields the weak form
\begin{eqnarray}
\rho \, c \, 
\left(\frac{\partial T}{\partial t}, \, \delta T \right)_{\Omega(t)} + \left( k \, \nabla T, \,
\nabla \delta T \right)_{\Omega(t)} &=&  \left(f , \,
 \delta T \right)_{\Omega(t)} 
+ \left(q_N , \,
 \delta T \right)_{\partial \Omega_N(t)} \nonumber \\ 
&+&   \left(  \sigma  , \,
\delta T \right)_{\Gamma(t)} +
\left(I \cdot  n_{\Gamma} , \,
\delta T \right)_{\Gamma(t)} 
\label{eq:stefan-weak signorini}
\end{eqnarray} 
with the modified Signorini conditions 
  \begin{equation}
  \begin{aligned}
\sigma  &\leq 0 &\mbox{ on } \Gamma(t), \\ 
T - T_m &\leq 0 &\mbox{ on } \Gamma(t), \\
\sigma \perp (T-&T_m)  & \mbox{ on } \Gamma(t).
  \end{aligned}
\label{equ: strong signorini conditions sigma}
  \end{equation}
We define bilinear form 
\begin{equation}
a(T,\delta T) : = \rho \, c \, 
\left(\frac{\partial T}{\partial t}, \, \delta T \right)_{\Omega(t)} + \left( k \, \nabla T, \,
\nabla \delta T \right)_{\Omega(t)}
\end{equation}
and linear form 
\begin{equation}
l(\delta T) :=  \left(f , \,
 \delta T \right)_{\Omega(t)} 
+ \left(q_N , \,
 \delta T \right)_{\partial \Omega_N(t)}
+ \left(I \cdot  n_{\Gamma} , \,
\delta T \right)_{\Gamma(t)}. 
\end{equation}
The weak formulation of the Stefan-Signorini problem then reads: For all $t \in [t_0,t_f]$, find $T \in H^1_D(\Omega(t))$ such that, for all $\delta T \in H^1_0(\Omega(t))$, 
\begin{equation}
a(T,\delta T) -  \left(  \sigma  , \,
\delta T \right)_{\Gamma(t)} = l(\delta T)
\label{equ: Signorini weak 1}
\end{equation}
with the Signorini law \eqref{equ: strong signorini conditions sigma} and \eqref{equ: latent heat interface} and initial conditions  \eqref{eq:initial condition}. Existence and uniqueness of the Stefan-Signorini problem are discussed in \cite{Friedman1984,Narimanyan2006}.

\subsection{Nonlinear Nitsche Reformulation}
\begin{figure}[htb]
  \begin{center}
    \includegraphics[width=0.5\textwidth]{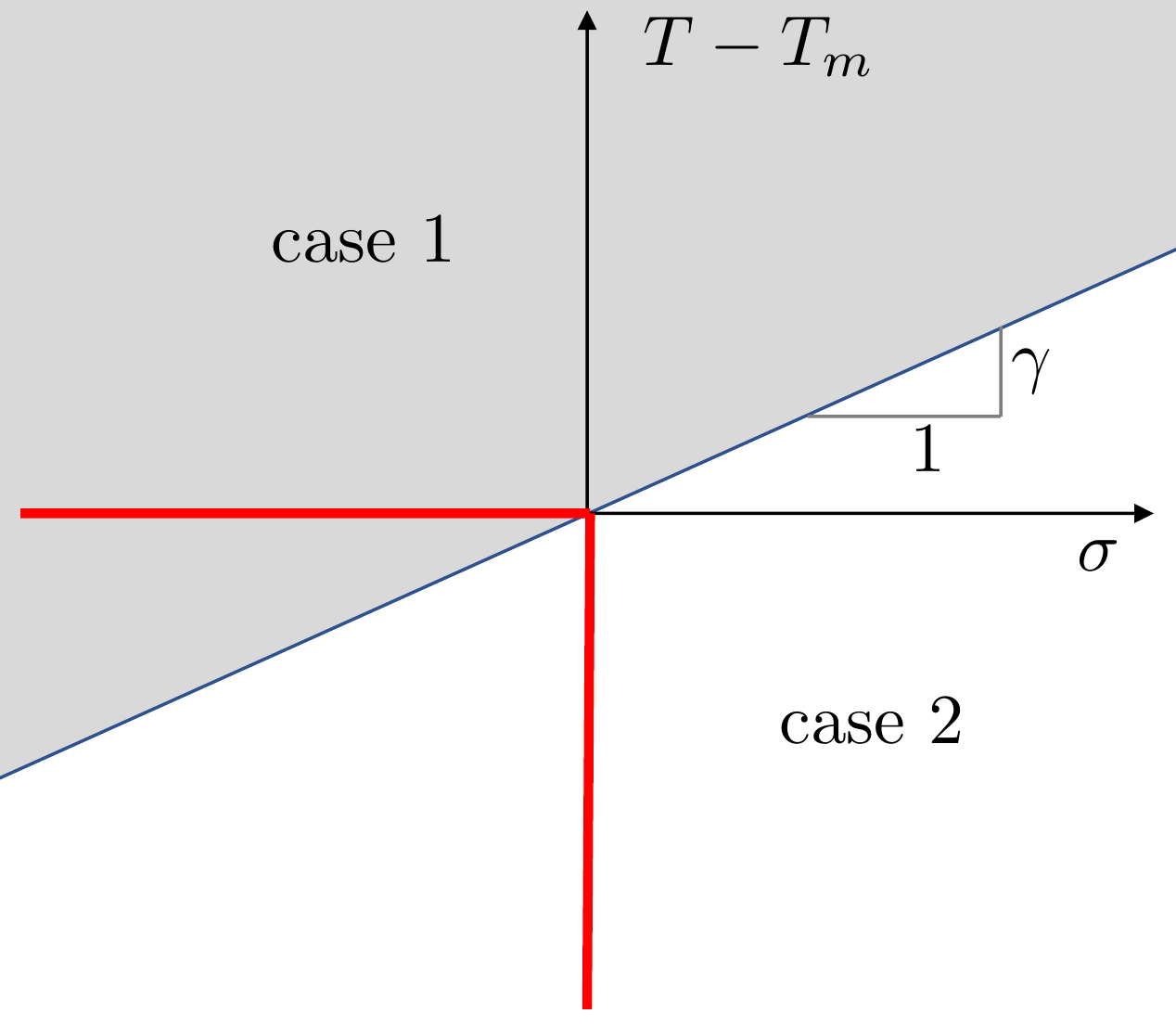}
  \end{center}
  \caption{Illustration of the different formulations of the Signorini law.}
\label{fig: Signorini illustration}
\end{figure}

First, following \cite{Curnier1988,Chouly2015,Burman2017a}, let us reformulate the Signorini law \eqref{equ: strong signorini conditions sigma} as
\begin{equation}
\sigma = - \frac{1}{\gamma} \left[ (T-T_m) - \gamma \sigma \right]_+ \, ,
\label{equ: Signorini relationship}
\end{equation}
where $\gamma$ is a positive penalty parameter and $[\cdot]_+$ denotes the positive part of a scalar quantity $x \in \mathbb{R}$, i.e. 
\begin{equation}
[x]_+ = 
\begin{cases}
x \quad \mbox{if } x>0, \\
0 \quad \mbox{otherwise}.
\end{cases}
\end{equation}
\begin{remark}
The equivalence of \eqref{equ: Signorini relationship} and \eqref{equ: strong signorini conditions sigma}  can be proved by enumeration. 

Let us first show that \eqref{equ: Signorini relationship} implies \eqref{equ: strong signorini conditions sigma}. Consider the following two complementary cases: \\
Case 1: $[(T-T_m) - \gamma \sigma] \geq 0$ $\Rightarrow$ $\sigma = - \frac{1}{\gamma} ( (T-T_m) - \gamma \sigma ) $ $\Rightarrow$ $T-T_m =0$. Now, returning to the first statement, $[0 - \gamma \sigma] \geq 0$ also indicates that $\sigma \leq 0$. As $\sigma \cdot (T-T_m ) = \sigma \cdot 0 = 0$, the Signorini law \eqref{equ: strong signorini conditions sigma} is satisfied in case 1.
\\
Case 2: $[(T-T_m) - \gamma \sigma] \leq 0$  $\Rightarrow$ $\sigma = 0$.  Now, returning to the previous statement, $[(T-T_m) - \gamma 0] \leq 0$ means that quantity $T-T_m$, which can be non-zero, is necessarily negative. Finally $\sigma \cdot (T-T_m ) =  0  \cdot (T-T_m )  = 0$, and therefore Signorini law \eqref{equ: strong signorini conditions sigma} is satisfied in case 2.
These two cases are illustrated in Figure~\ref{fig: Signorini illustration}. 

Let us now show that \eqref{equ: strong signorini conditions sigma} implies \eqref{equ: Signorini relationship}. Consider the three following cases: \\
Case 1: $\sigma < 0$. Owing to the consistency condition, this can only happen when $T-T_m = 0$. Therefore, we can write $\sigma = -\frac{1}{\gamma}(-\gamma \sigma) = -\frac{1}{\gamma}[0-\gamma \sigma]_+ =  -\frac{1}{\gamma}[T-T_m-\gamma \sigma]_+$. 
\\
Case 2: $T-T_m < 0$. This can only happen when $\sigma=0$. Therefore, we can write that $\sigma = 0 =  -\frac{1}{\gamma}[T-T_m]_+=  -\frac{1}{\gamma}[T-T_m-\gamma \sigma]_+$.
\\
Case 3: The last possible scenario is $\sigma = T-T_m =0$. In this case, we can write that $\sigma = 0 = -\frac{1}{\gamma}[ 0 ]_+ =  -\frac{1}{\gamma}[T-T_m-\gamma \sigma]_+$. $\hfill \square$

\end{remark}

A Nitsche formulation of the Stefan-Signorini problem~\eqref{equ: Signorini weak 1} is obtained by replacing slack variable $\sigma$ by its expression as a function of the primal variable $T$,
\begin{equation}
\sigma(T) =  - \frac{1}{\gamma} \left[ (T-T_m) - \gamma (k \nabla T \cdot n_{\Gamma} - I \cdot n_{\Gamma}) \right]_+ \, .
\end{equation}
A penalty term enforcing this expression weakly can be formulated as follows
\begin{equation}
 s_\heartsuit(T,\delta T)  := \int_{\Gamma(t)}   \left( ( k \, \nabla T \cdot n_{\Gamma} - I \cdot  n_{\Gamma}) - \sigma(T) \right) \left(  \theta_1 \delta T - \theta_2 \gamma k \nabla \delta T \cdot n_{\Gamma} \right)  d \Gamma\, ,
\end{equation}
where the choice of $\theta_1 \in \{0,1\}$ and $\theta_2 \in \{ -1, 0, 1\}$ leads to a family of different methods as detailed below. Finally, the proposed \textit{Stefan-Signorini-Nitsche} formulation reads 
\begin{equation}
a(T,\delta T) -  \left( k \, \nabla T \cdot n_{\Gamma} - I \cdot  n_{\Gamma}  , \,
\delta T \right)_{\Gamma(t)} + s_\heartsuit(T,\delta T) = l(\delta T) \,.
\end{equation}
To simplify the notations, let us define
\begin{equation}
P_{\gamma}(T) := (T-T_m) - \gamma (k \nabla T \cdot n_{\Gamma} - I \cdot n_{\Gamma}) \, ,
\end{equation}
and the parametrised variation of this quantity, which we define as
\begin{equation}
P^\delta_{\theta\gamma}(\delta T) := \theta_1 \delta T -  \gamma \, \theta_2 \, k \, \nabla \delta T \cdot n_{\Gamma} \, .
\end{equation}
Using these notations, the penalty term $s_\heartsuit$ reads as 
\begin{equation}
 s_\heartsuit(T,\delta T)  = \int_{\Gamma(t)}   \left( ( k \, \nabla T \cdot n_{\Gamma} - I \cdot  n_{\Gamma}) + \frac{1}{\gamma} \left[ P_{\gamma}(T) \right]_+  \right)  P^ \delta_{\theta\gamma}(\delta T)   d \Gamma\, .
\end{equation}
And the proposed Nitsche formulation can be expressed as a sum of linear and nonlinear terms, 
\begin{equation}
 a(T,\delta T)  
+
 \left( k \, \nabla T \cdot n_{\Gamma},  P^\delta_{\theta\gamma}(\delta T) - \delta T \right)_{\Gamma(t)} + \mathcal{N}(T,\delta T) 
= l(\delta T) 
+
 \left( I \cdot n_{\Gamma}, P^\delta_{\theta\gamma}(\delta T)  - \delta T 
\right)_{\Gamma(t)} \, ,
\label{equ: reformulated Nitsche Stefan Signorini}
\end{equation}
where
\begin{equation}
\mathcal{N}(T,\delta T) =   \frac{1}{\gamma} \left( \left[ P_{\gamma}(T) \right]_+, P^\delta_{\theta\gamma}(\delta T) 
\right)_{\Gamma(t)}  \, .
\end{equation}
We emphasise the fact that $\mathcal{N}$ is nonlinear in its first argument.

Three interesting Nitsche formulations are obtained by choosing particular values for the $(\theta_1,\theta_2)$  pair, as follows.
\begin{itemize}
\item{$\theta_1 = 1$ and $\theta_2 = 1$, $\gamma >0$:} Symmetric Nitsche method.
\item{$\theta_1 = 1$ and $\theta_2 = -1$, $\gamma >0$:} Non-symmetric Nitsche method. A closely related formulation was proposed in \cite{Chouly2015} to solve problems of unilateral contact between deformable solids. It was shown that, as opposed to the symmetric Nitsche formulation, the stability of this non-symmetric variant is preserved irrespectively of the value of Nitsche parameter $\gamma$. 
\item{$\theta_1 = 1$ and $\theta_2 = 0$, $\gamma >0$:} Consistent penalty formulation as described in \cite{Curnier1988}.
\item{$\theta_1 = 0$ and $\theta_2 = -1$:}  Semi penalty-free Nitsche method. A closely related formulation was derived in the context of Signorini-Poisson problems in \cite{Burman2017a}. In this setting, the method is stable irrespectively of $\gamma$. The term ``penalty-free" indicates that condition $T=T_m$ is enforced with the penalty-free Nitsche method (i.e. non-symmetric Nitsche method without penalty term). However, $\gamma$  is the scaling of a penalty term that enforces the Neumann interface condition when $T<T_m$.
\end{itemize}
In this article, we focus on the (semi) penalty-free Nitsche method. Our motivation is that, at least in the context of equality constraints, the penalty-free Nitsche method yields better interface fluxes than the symmetric and non-symmetric Nitsche method, as was shown in \cite{Boiveau2017}. Interface fluxes are of particular importance in the Stefan-Signorini problem because they drive the motion of the interface when $T=T_m$.

\section{Stabilised Cut Finite Element Formulation}
\label{sec:fem-formulation}
In this section, we introduce the spatial and temporal finite element discretisation of problem~\eqref{equ: reformulated Nitsche Stefan Signorini}.
\subsection{Discretisation in Space}
\subsubsection{Background Mesh and Fictitious Domain}
\begin{figure}[t!]
  \begin{center}
    \includegraphics[width=0.7\textwidth]{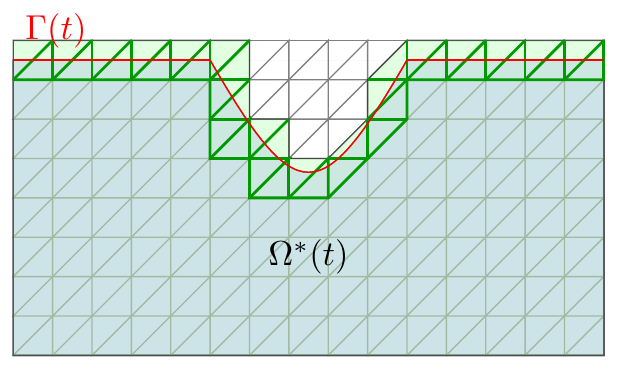}
  \end{center}
  \caption{Schematics of  the domain $\Omega(t)$ covered by a fixed and regular background mesh $\tilde{\mesh}$ and the fictitious domain $\Omega^*(t)$ consisting of all elements in $\tilde{\mesh}$ with at least one part in $\Omega(t)$. 
}
  \label{fig:computational-domain}
\end{figure}

First, let us introduce important background mesh quantities and the definition of the evolving fictitious domain. Let $\Omega(t_0)$ be our domain in $\R^d$ ($d = 2, 3$) at time $t=t_0$ with
Lipschitz boundary $\partial \Omega$ and let $\tilde{\mesh}$ be a
quasi-uniform tesselation that covers the domain $\Omega(t_0)$.
We define a background domain 
\begin{equation}
\Omega_b = \bigcup_{K \in \tilde{\mesh}} K
\end{equation}
associated with our fixed tesselation $\tilde{\mesh}$. The background mesh $\tilde{\mesh}$ will stay fixed in time while the interface $\Gamma(t)$ will move through this background mesh. For $t \in [t_0,t_f]$, we denote the elements in the background mesh $\tilde{\mesh}$ that have at least a small part in domain $\Omega(t)$ as
\begin{equation}
\mesh(t) = \{ K \in \tilde{\mesh} : K \cap \Omega(t) \neq \emptyset\},
\end{equation}
which we call \emph{active mesh} (see the blue and green shaded elements in Figure~\ref{fig:computational-domain}). In contrast to the background mesh, the active mesh changes in time. We denote the union of all elements in $\mesh(t)$ as 
\begin{equation}
\Omega^*(t) = \bigcup_{K \in \mesh(t)} K.
\end{equation}
$\Omega^*(t)$ is called the \emph{fictitious domain}. \\ 
In addition to these time-evolving domains and meshes, we have edge stabilisation quantities that change in time. For each active mesh at time $t\in [t_0,t_f] $, we will distinguish between the following different sets of faces, i.e. edges in 2D and faces in 3D. The \emph{exterior faces}, $\mathcal{F}_e(t)$, which are the faces that belong to one element only in the background mesh and that have an intersection with the active mesh. The \emph{interior faces}, $\mathcal{F}_i(t)$, which are faces that are shared by two elements with $K \cap \Omega(t) \neq \emptyset$. 
To prevent ill-conditioning, we will apply stabilisation terms to the elements which are intersected by the boundary $\Gamma(t)$, i.e.
\begin{equation}
  \label{eq:define-cutting-cell-mesh}
G_h(t) = \{K \in \mesh(t): K \cap \Gamma(t) \neq \emptyset \}.
\end{equation}
These stabilisation terms will be applied to so-called \emph{ghost penalty} faces defined as 
\begin{equation}
  \label{eq:define-cutting-face-mesh}
  \mcF_{\Gamma}(t) = \{ F \in \mcF_i(t) :\;
  K^+_F \cap \Gamma(t) \neq  \emptyset
  \vee
  K^-_F \cap \Gamma(t) \neq  \emptyset
  \}.
\end{equation}
Here, $K^+_F$ and $K^-_F$ are the two elements sharing the interior
face $F \in \mathcal{F}_i(t)$. The set of faces $ \mcF_{\Gamma}(t)$ is illustrated by the dark green edges shown in Figure~\ref{fig:computational-domain}. To ensure that the boundary $\Gamma(t)$ is reasonably resolved by $\mesh$, we make the following assumptions:
\begin{itemize}
\item G1: The intersection between $\Gamma(t)$ and a face $F \in
  \mcF_i(t)$ is simply connected; that is, $\Gamma(t)$ does not
  cross an interior face multiple times.
\item G2: For each element $K$ intersected by $\Gamma(t)$, there exists a
  plane $S_K$ and a piecewise smooth parametrization $\Phi: S_K \cap K
  \rightarrow \Gamma(t) \cap K$.
\item G3: We assume that there is an integer $N>0$ such that for each
  element $K \in G_h(t)$ there exists an element $K' \in
  \mesh(t) \setminus G_h(t)$ and at most $N$ elements
  $\{K\}_{i=1}^N$ such that $K_1 = K,\,K_N = K'$ and $K_i \cap K_{i+1}
  \in \mcF_i(t),\; i = 1,\ldots N-1$.  In other words, the
  number of faces to be crossed in order to ``walk'' from a cut
  element $K$ to a non-cut element $K' \subset \Omega(t)$ is uniformly bounded.
\end{itemize}
Similar assumptions were made in \citep{HansboHansbo2002,BurmanHansbo2012}.

\subsubsection{Nonconforming spatial discretisation of the Stefan-Signorini problem}
Using the sets of mesh elements and faces defined above, we can formulate the discrete Stefan-Signorini problem. Firstly, we introduce the continuous linear finite element
space on the \emph{active} mesh 
\begin{align}
  \mathcal{V}_h(t) &= \left\{ v_h \in C^0(\Omega^*(t)): \restr{v_h}{K} \in
    \mathcal{P}_1(K)\, \foralls K \in \mesh(t) \right\}
\end{align}
for the temperature. \\
Secondly, we define a stabilisation operator on the faces $ \mcF_{\Gamma}(t)$ for the temperature $T$ to prevent ill-conditioning in the case of intersections of $\Gamma(t)$ near a node or face of elements as
\begin{equation}
  \label{eq:stabilization_operator}
  s_{T}(T_h,\delta T_h) = \sum\limits_{F \in \mathcal{F}_{\Gamma}(t)} \gamma_T \, k \, h  
  \left( \jump{\nabla T_h}_n, \jump{\nabla \delta T_h}_n \right)_F.
\end{equation} 
Here, $\jump{\nabla x}_n$ denotes the normal jump of the quantity $x$
over the face, $F$, defined as $\jump{\nabla x}_n = \left. \nabla x
\right|_{T_F^+} n_F  - \left. \nabla x
\right|_{T_F^-} n_F$, where $n_F$ denotes a unit
normal to the face $F$ with fixed but arbitrary orientation and 
 $\gamma_T$ is a positive penalty parameter
to be determined later. We refer to the term $ s_{T}(T_h,\delta T_h)$ as \textit{ghost penalty} stabilisation \cite{Burman2010}. 
Using the definitions above, we are now in the position to formulate our stabilised cut finite
element method for the one-phase Stefan-Signorini problem.
The proposed discretisation scheme reads: For all $t \in [t_0,t_f]$, find
$T_h \in   \mathcal{V}_h(t)$ such that for all
$\delta T_h  \in   \mathcal{V}_h(t)$
\begin{equation}
\begin{aligned}
A(T_h,\delta T_h) + \mathcal{N}(T_h,\delta T_h) &= L(\delta T_h) \, , 
\label{equ: discretised weak form}
\end{aligned}
\end{equation}
where
\begin{equation}
\begin{aligned}
A(T_h,\delta T_h)  &= a(T_h,\delta T_h)  +     a_{bc}(T_h, \delta T_h) + s_{T}(T_h,\delta T_h) +
 \left( k \, \nabla T_h \cdot n_{\Gamma},  P^\delta_{\theta\gamma}(\delta T_h) - \delta T_h \right)_{\Gamma_h(t)} \, , \\
L(\delta T_h) &= l(\delta T_h) +   l_{bc}(T_h, \delta T_h)
+
 \left( I \cdot n_{\Gamma}, P^\delta_{\theta\gamma}(\delta T_h)  - \delta T_h 
\right)_{\Gamma_h(t)} \, , \\
 \mathcal{N}(T_h,\delta T_h) &=   \frac{1}{\gamma} \left( \left[ P_{\gamma}(T_h) \right]_+, P^\delta_{\theta\gamma}(\delta T_h) 
\right)_{\Gamma_h(t)}  \, , \\
\end{aligned}
\end{equation}
with 
  \begin{align}
  \label{eq:A_h-definition}
    a(T_h, \delta T_h) =& \rho \, c \, 
\left(\frac{\partial T_h}{\partial t}, \, \delta T_h \right)_{\Omega_h(t)} + \left( k \, \nabla T_h, \,
\nabla \delta T_h \right)_{\Omega_h(t)},   \\
  \label{eq:L_h-definition}
  l(\delta T_h) =&  \left(f , \,
 \delta T_h \right)_{\Omega_h(t)} 
+ \left(q_N , \,
 \delta T_h \right)_{\partial \Omega_N(t)}
+ \left(I \cdot  n_{\Gamma} , \,
\delta T_h \right)_{\Gamma_h(t)}  , \\
    a_{bc}(T_h, \delta T_h) =& - \left( k \, \nabla T_h \cdot n_{\Omega}, \,
 \delta T_h \right)_{\partial \Omega_D(t)}  - \left( k \, \nabla \delta T_h \cdot n_{\Omega}, \,
 T_h \right)_{\partial \Omega_D(t)}  + \frac{k \gamma_b}{h} \left(  \, T_h , \,
 \delta T_h \right)_{\partial \Omega_D(t)}, \\ 
    l_{bc}(T_h, \delta T_h) =&   - \left( k \, \nabla \delta T_h \cdot n_{\Omega}, \,
 T_D \right)_{\partial \Omega_D(t)}  + \frac{k \gamma_b}{h} \left(  \, T_D , \,
 \delta T_h \right)_{\partial \Omega_D(t)}. 
  \end{align}
Here, the positive penalty constant $\gamma_b$ arises from the weak
enforcement of Dirichlet boundary conditions through Nitsche's method
and $h=\max_{K\in \mesh} h_K$ is the maximum mesh size, where $h_K$ denotes
the diameter of $K$. The penalty parameter $\gamma$ is now scaling as $\gamma = \hat{\gamma} h$ with $\hat{\gamma}>0$ chosen to be sufficiently small. 

The discretised Stefan-Signorini problem is completed with initial conditions 
  \begin{equation}
\label{eq:InitialCondInterp}
  \begin{aligned}
T_h(x,t_0) &=  \hat{\mathcal{I}} (T_0) &\mbox{ in } \Omega_h(t_0),
  \end{aligned}
  \end{equation}
and the condition on the normal velocity of the boundary $\Gamma(t)$
  \begin{equation}
v \cdot n_{\Gamma}  = \frac{k \, \nabla T_h \cdot n_{\Gamma} - I  \cdot n_{\Gamma}}{\rho \, L},
\label{equ: normal velocity interface 1}
  \end{equation}
whose discretisation will be discussed in detail in Section 5. In equation \eqref{eq:InitialCondInterp}, $\hat{\mathcal{I}}$ is the standard finite element nodal interpolation operator. Note that it is implicitly assumed that field $T_0$ is analytically available in the entire fictitious domain.

\subsection{Discretisation in Time} 
We decompose the time interval $[t_0,t_f]$ into $n_t$ time steps, and 
we seek a sequence of solutions $\{ T( t_n)\}_{n \in \llbracket 0 \ n_t -1\rrbracket} =:\{ T^n \}_{n \in \llbracket 0 \ n_t -1 \rrbracket} $. The reference time $t_0$ is chosen to be equal to 0. 
We assume that times $\{ t_n \}_{n \in \{ 0, \dots, \ n_t-1 \} } $ are uniformly spaced, which allows us to define the time step $\Delta t = t_1 - t_0 (=t_2-t_1=...)$. We apply a backward Euler scheme to the system \eqref{equ: discretised weak form} and evaluate integrals over the domain $\Omega_h(t_n)$ and the boundary $\Gamma_h(t_n)$. This time discretisation yields the fully discrete system of equations at time step $n+1$: Find $T^{n+1} \in \mathcal{V}_h(t_n)$, such that for all  $\delta T \in \mathcal{V}_h (t_n)$ 
\begin{equation}
A_{\sharp}(T_h^{n+1},\delta T) + \mathcal{N}(T_h^{n+1},\delta T) = L_{\sharp}(\delta T) \, ,
\label{equ: fully discretised system}
\end{equation}
where
\begin{equation}
\begin{aligned}
A_{\sharp}(T_h^{n+1},\delta T_h)  &= a_{\sharp}(T_h^{n+1},\delta T_h)  +     a_{bc}(T_h^{n+1}, \delta T_h) \\
+& s_{T}(T_h^{n+1},\delta T_h) +
 \left( k \, \nabla T_h^{n+1} \cdot n_{\Gamma},  P^\delta_{\theta\gamma}(\delta T_h) - \delta T_h \right)_{\Gamma_h(t_n)} \, , \\
  L_{\sharp}(\delta T_h) &=  L(\delta T) +   \rho \, c \, 
\left(\frac{T_h^{n}}{\Delta t}, \, \delta T_h \right)_{\Omega_h(t_n)}  
\end{aligned}
\end{equation}
with 
  \begin{align}
  \label{eq:A_h-definition}
    a_{\sharp}(T_h^{n+1}, \delta T_h) =& \rho \, c \, 
\left(\frac{ T_h^{n+1}}{\Delta t}, \, \delta T_h \right)_{\Omega_h(t_n)} + \left( k \, \nabla T_h, \,
\nabla \delta T_h \right)_{\Omega_h(t_n)} . 
  \end{align}

\subsection{Newton-Raphson algorithm} 

We solve the previous system \eqref{equ: fully discretised system} using a semi-smooth Newton-Raphson algorithm. We linearise the  semi-linear form $ \mathcal{N}$ around a finite element reference temperature $T_h^\star \in \mathcal{V}_h(t_n)$. This is done by writing the Taylor expansion, for any $\delta T \in \mathcal{V}_h(t_n)$,
\begin{equation}
 \mathcal{N}(T_h,\delta T) =  \mathcal{N}(T_h^\star,\delta T) +  D\mathcal{N}(T_h - T_h^\star ,\delta T;T_h^\star) \, ,
\end{equation}
where $T_h \in \mathcal{V}_h(t_n)$ and  
$\mathcal{D}\mathcal{N}$ is the G\^{a}teaux-derivative of $\mathcal{N}$, which is defined, for any finite element field $d{T}_h \in \mathcal{V}_h(t_n)$, by
\begin{equation}
\mathcal{D}\mathcal{N}( d{T}_h,\delta T;T_h^\star) = \lim_{z \rightarrow 0} \frac{1}{z} \left(  \mathcal{N}(T_h^\star+z \, d T_h,\delta T) -  \mathcal{N}(T_h^\star,\delta T) \right) \, .
\end{equation}
Identifying the temperature increment $d{T}_h = T_h - T_h^\star \in \mathcal{V}_h(t_n)$, we find that
\begin{equation}
\mathcal{D}\mathcal{N}(d{T}_h,\delta T;T_h^\star) = 
\frac{1}{\gamma} \left( \mathcal{D}G( d{T}_h; T_h^\star) , P^\delta_{\theta\gamma}(\delta T) 
\right)_{\Gamma(t_n)}  \, ,
\end{equation}
where $\mathcal{D}G( \  . \ ; T_h^\star)$ is the G\^{a}teaux-derivative of $G(T):= \left[ P_{\gamma}(T) \right]_+$ at $T_h^\star$, which is given by
\begin{equation}
\mathcal{D}G(d{T}_h; T_h^\star) = H(P_{\gamma}(T_h^\star))  \left(d{T}_h - \gamma k \nabla d{T}_h \cdot n_\Gamma \right) \, ,
\end{equation}
where $H$ is the Heaviside function
\begin{equation}
H(P_{\gamma}(T_h^\star)) = 
\begin{cases}
1 \quad \mbox{if } P_{\gamma}(T_h^\star) >0, \\
0 \quad \mbox{otherwise}. 
\end{cases}
\end{equation}
Using these derivations, the Newton predictor for the $(k+1)^\text{th}$ iterate $T_h^{k+1}$ of $T_h^{n+1}$ is, for any $k \in \mathbb{N}^+$ and for all $\, \delta T \in \mathcal{V}_h(t_n)$, given by
\begin{equation}
A_{\sharp}( d{T}_h,\delta T) + \mathcal{D}\mathcal{N}(d{T}_h,\delta T;T_h^k) 
= r(\delta T;T_h^k) \, ,
\end{equation}
where the Newton increment is defined as $d{T}_h:=T_h^{k+1} - T_h^{k}$ and the residual $r$ of iterate $T_h^k$ is such that for any $\delta T \in \mathcal{V}_h(t_n)$
\begin{equation}
 r(\delta T;T_h^k) = L_{\sharp}(\delta T) - \left(  A_{\sharp}(T_h^{k},\delta T) + \mathcal{N}(T_h^k,\delta T) \right).
\end{equation}

\section{Description of the Domain Movement}
In this Section, we describe how domain $\Omega(t)$ is discretised and evolved in time. For each time-step, $t \in [t_0,t_f]$, we first solve  \eqref{equ: fully discretised system} to obtain the temperature $T_h^{n+1}$, with which we determine the normal velocity on $\Gamma(t_n)$, i.e. 
\begin{equation}
v^{n+1} \cdot n_{\Gamma}  = \frac{k \, \nabla T^{n+1} \cdot n_{\Gamma} - I  \cdot n_{\Gamma}}{\rho \, L}  \quad \mbox{on } \Gamma(t_n).
\label{equ: normal velocity interface}
\end{equation}
Then, this normal velocity on $\Gamma(t_n)$ is used to move a level-set function as detailed below. 
\subsection{Level-set description of the moving domain and level-set advection}
We track the motion of the boundary $\Gamma(t)$ using a continuous level-set function $\phi: \Omega_b \times [t_0,t_f] \rightarrow \mathbb{R}$, whose zero level set describes the location of the  boundary $\Gamma(t) = \{x \in \Omega_b: \phi(x,t)=0,t \in [t_0,t_f] \}$. The material domain $\Omega(t)$ is implicitly defined by $\phi(x,t)<0$, and the fictitious domain by $\phi(x,t)>0$. \\
To satisfy equation~\eqref{equ: normal velocity interface}, the zero level-set contour is required to move with $v^{n+1} \cdot n_{\Gamma}$. 
Furthermore, as the level-set function is defined over the entire background domain $\Omega_b$, the velocity at the interface \eqref{equ: normal velocity interface} needs to be extended to the remainder of the domain (or at least to a band around the zero isoline of the level-set) to evolve the level-set function. We denote this extension of the velocity with $v_{ext}$. \\
Then, the level-set function is moved by solving the advection problem 
\begin{equation}
\frac{\partial \phi}{\partial t} + v_{ext} \cdot \nabla \phi = 0
\label{equ: level set advection strong}
\end{equation}
with initial condition $\phi(x,0)=\phi_0$. Here, $\phi_0$ is a given initial level-set description of the domain $\Omega(t_0)$.
We discretise the level set function using a continuous quadratic finite element space defined on the entire background mesh $\tilde{\mesh}$ denoted by
\begin{equation}
\mathcal{W}_h := \left\{ v_h \in C^0(\Omega_b): \restr{v_h}{K} \in
    \mathcal{P}_2(K)\, \foralls K \in \tilde{\mesh} \right\}.
\end{equation}
To solve the advection equation~\eqref{equ: level set advection strong}, we use a $\theta$-scheme in time and streamline diffusion (SUPG) in space. The discretised advection problem reads: Find $\phi_h^{n+1} \in \mathcal{W}_h$, such that for all $\delta \phi \in \mathcal{W}_h$
\begin{equation}
a_{\phi}(\phi_h^{n+1},\delta \phi) = l_{\phi}(\delta \phi)
\label{equ: level set advection}
\end{equation}
with 
\begin{equation}
\begin{aligned}
a_{\phi}(\phi_h^{n+1},\delta \phi) &= \left( \frac{\phi_h^{n+1}}{\Delta t}+ \theta v_{ext}^{n+1} \cdot \nabla \phi_h^{n+1}, \delta \phi + \tau_{SD} (v_{ext} ^{n+1} \cdot \nabla \delta \phi) \right)_{\Omega_b}, \\
l_{\phi}(\delta \phi) &=  \left( \frac{\phi_h^{n}}{\Delta t}+ (1-\theta) v_{ext}^{n} \cdot \nabla \phi_h^{n}, \delta \phi + \tau_{SD} (v_{ext} ^{n+1} \cdot \nabla \delta \phi) \right)_{\Omega_b}
\label{equ: level set SUPG}
\end{aligned}
\end{equation}
with the streamline diffusion parameter (see \cite{hansbo2016})
\begin{equation}
\tau_{SD} = 2 \left( \frac{1}{\Delta t^2} + \frac{v_{ext} \cdot v_{ext} }{h^2} \right)^{-\frac{1}{2}}
\end{equation}
and initial condition $\phi_h^0 =\phi_0$.
Throughout this contribution, $\theta$ is set to 0.5. 
\subsection{Description of the geometry}
\begin{figure}
\centering
\includegraphics[width=.5\textwidth]{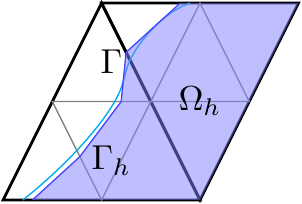}
\caption{Illustration of linear approximation of interface $\Gamma$ and domain $\Omega$ on the refined mesh $\tilde{\mesh}_{/2}$ with respect to the mesh $\tilde{\mesh}$.}
\label{fig: refined mesh approximation}
\end{figure}

The quadratic level-set function is used to define the geometry of our problem including the discretised material domain  $\Omega_h(t)$, the discretised interface $\Gamma_h(t)$ and the normal $n_{\Gamma}(t)$. In the rest of this section we choose a fixed time t and suppress the time dependence to ease the notation.
\paragraph{Normal Computation} 
The normal pointing from the domain $\Omega(t)$ into the fictitious material at the interface $\Gamma(t)$ can be obtained from the level-set function using 
\begin{equation}
n_{\Gamma}(x,t) = \frac{\nabla \phi(x,t)}{\norm{\nabla \phi(x,t)}}. 
\end{equation}
In this contribution, we determine the normal $n_{\Gamma}$ from the level set function through a $L^2$-projection onto the continuous piecewise linear space 
\begin{equation}
\mathcal{X}_h^d := \left\{ v_h \in [C^0(\Omega_b)]^d : \restr{v_h}{K} \in \mathcal{P}_1(K)\, \forall K \in \tilde{\mesh} \right\}.
\end{equation} 
Here, $d=2,3$ is the geometrical dimension. We determine the normal by finding $n_{\Gamma} \in \mathcal{X}_h^d $ such that for all $\delta n_{\Gamma}  \in \mathcal{X}_h^d$
\begin{equation}
(n_{\Gamma},\delta  n_{\Gamma})_{\Omega_b} = \left(\frac{\nabla \phi}{|\nabla \phi|},\delta  n_{\Gamma}\right)_{\Omega_b}. 
\label{equ: l2 projection normal}
\end{equation} 
\paragraph{Discrete Geometrical Domains}
To define our discrete domains $\Omega_h$ and $\Gamma_h$, we use a two-grid solution proposed by \cite{Gross2006, Gross2011} which is outlined in the following. First, we interpolate the piecewise quadratic level-set function onto a piecewise linear function, $\mathcal{I}(\phi_h)$ , on a regularly refined mesh $\tilde{\mesh}_{/2}$, such that  
\begin{equation}
\mathcal{I}(\phi_h(\mathbf{v})) = \phi_h(\mathbf{v}) \mbox{ for all nodes $\mathbf{v}$ in }\tilde{\mesh}_{/2}.
\end{equation}
We then use this piecewise linear interpolation to determine the intersection between $\mathcal{I}(\phi_h)$ with the refined grid to obtain the piecewise linear approximation of $\Omega_h$ and $\Gamma_h$ as illustrated in Figure~\ref{fig: refined mesh approximation}.
\subsection{Interface velocity smoothing and bulk extension}
To enable a smooth domain movement, we construct a continuous piecewise linear normal velocity approximation in the vicinity of the interface in the following way. We first recover a smoothed gradient of the temperature using the stabilised projection: Find $G_T^h \in  \mathcal{V}_h^d(t)$ such that for all $\delta G_T \in \mathcal{V}_h^d(t)$
\begin{align}
\left( G_T^h, \delta G_T \right)_{\Omega_h} + s_{G_T}(G_T^h, \delta G_T) &= (\nabla T^{n+1}_h, \delta G_T)_{\Omega_h}, \\ 
 s_{G_T}(G_T^h, \delta G_T) &= \sum\limits_{F \in \mathcal{F}_{\Gamma}(t)} \gamma_{G_T} \, \, h   \left( \jump{\nabla {G_{T}^h}}_n, \jump{\nabla \delta G_T}_n \right)_F,
\label{equ: smoothed gradient projection}
\end{align}
where 
\begin{align}
  \mathcal{V}_h^d(t) &= \left\{ v_h \in [C^0(\Omega^*(t))]^d: \restr{v_h}{K} \in
    \mathcal{P}_1(K)\, \foralls K \in \mesh(t) \right\}
\end{align}
is the vector-valued space of continuous piecewise linear functions on the fictitious domain with $d=2,3$. Here, $\gamma_{G_T} >0$ is a positive penalty parameter to recover a continuous piecewise linear temperature gradient over patches of elements in the interface region from the piecewise constant temperature gradient $\nabla T^{n+1}_h$. The normal velocity is then determined from the following $L^2$ projection: Find $v_n^h := v^{n+1}_h \cdot n_{\Gamma} \in \mathcal{V}_h(t)$ such that for all $\delta v_n \in  \mathcal{V}_h(t)$
\begin{equation}
\left( v_n^h, \delta v_n \right)_{\Omega^*}  = H(P_{\gamma}(T_h^{n+1})) \left[ \left(\frac{(k G_T^h-I) \cdot n_{\Gamma}}{\rho L}, \delta v_n\right)_{\Omega^*} - \frac{\theta_1}{\gamma \rho L}(T_h^{n+1}-T_m,\delta v_n)_{\Gamma_h} \right].
\label{equ: normal vel computation}
\end{equation}
We extend this normal velocity field onto the entire background domain $\Omega_b$ using a fast marching scheme as detailed in \cite{Jahn2017,Jahn2016,Gross2011}. The principle of this technique relies on two steps: a near-field step and a far-field  step. In the near-field step, all nodes of elements which are intersected by the interface obtain the normal velocity value determined by a closest point projection onto the discretised interface $\Gamma_h$, i.e. $v_n^{ext}(\mathbf{v}) = v_n(P_W(\mathbf{v}))$, where $\mathbf{v}$ is an element node and $P_W$ is the closest point projection onto $p_{\Gamma}(\mathbf{v}) \in \Gamma_h$ given by the shortest distance between $\mathbf{v}$ and $\Gamma_h$. 
In the far-field step the information of the intersected elements is propagated through evaluation of known extended velocity function values. For a detailed description see \cite{Jahn2017,Jahn2016,Gross2011}. \\
This extended normal velocity field, $v_n^{ext}$, is then used to obtain the vectorial extended velocity field as 
\begin{equation}
v_{ext} = v_n^{ext} \cdot n_{\Gamma},
\end{equation}
which is then used to propagate the level-set function using equation~\eqref{equ: level set SUPG}. 
\begin{remark}
This extended  velocity field enables the transport of the level-set function in a way that prevents a large deviation of the level set away from a signed distance function. Therefore reinitialisation of the level set is rarely required. However, in examples with strong interface deformation reinitialisation may become necessary. In these rare cases, we use a fast marching redistancing technique described in \cite{Gross2011,Gross2006}, which relies on a fast marching scheme of the linearly interpolated level set function on the regularly refined grid $\tilde{\mesh}_{/2}$.
\end{remark}

\noindent Algorithm~\ref{alg: cutfem stefan signorini} summaries the CutFEM Stefan-Signorini algorithm presented in the previous sections. 

\begin{algorithm}
    \caption{CutFEM Stefan-Signorini Algorithm}
    \label{alg: cutfem stefan signorini}
    \begin{algorithmic}[1] 
    \State Set $t=t_0$, $T_h^0=T_0$, $\phi_h^0 = \phi_0$.
    \While{$t\leq t_f$} \Comment{$t_f$ is the final time}
   \State Determine $\Omega_h$ and $\Gamma_h$ through intersection computations of zero level-set with background mesh.   
\State Compute normal $n_{\Gamma}$ using \eqref{equ: l2 projection normal}.
  \Procedure{Stefan-Signorini-Nitsche}{$T_h^n$} 
	\State Solve \eqref{equ: fully discretised system} using Newton-Raphson algorithm.
	 \State \textbf{return}	$T^{n+1}_h$. 
  \EndProcedure
  \Procedure{Velocity}{$T_h^{n+1}$}
	\State Compute smoothed temperature gradient $G_T$ using \eqref{equ: smoothed gradient projection}.
	\State Determine normal velocity on $\Gamma_h$ using \eqref{equ: normal vel computation}.
	\State Extend normal velocity to obtain $v_{ext}$. 
	\State \textbf{return}	$v_{ext}$. 
\EndProcedure
 \Procedure{Level Set Advection}{$v_{ext}$}
\State Solve advection problem~\eqref{equ: level set advection} for level set. 
\State \textbf{return}	$\phi_{h}^{n+1}$. 
\EndProcedure
       \State $t =t+\Delta t$, $\phi_h^{n} = \phi_{h}^{n+1}$, $T^n_h = T^{n+1}_h$.
    \EndWhile\label{while}
    \end{algorithmic}
\end{algorithm}

\section{Numerical results}

In this Section, we present numerical results for a manufactured solution, and for several thermal ablation problems in 2D and 3D.  The penalty parameters are set to $\gamma_{G_T}=10^{-3}$, $\gamma_{T}=10^{-1}$, $\hat{\gamma}=1$, $\gamma_b=100$ and we choose the penalty-free, non-symmetric Nitsche method, i.e. $\theta_1=0$, $\theta_2=-1$ in all the presented results.

\subsection{Manufactured solution: convergence analysis}
\begin{figure}
\begin{center}
\subfloat[Schematics.]{\includegraphics[width=.59\textwidth]{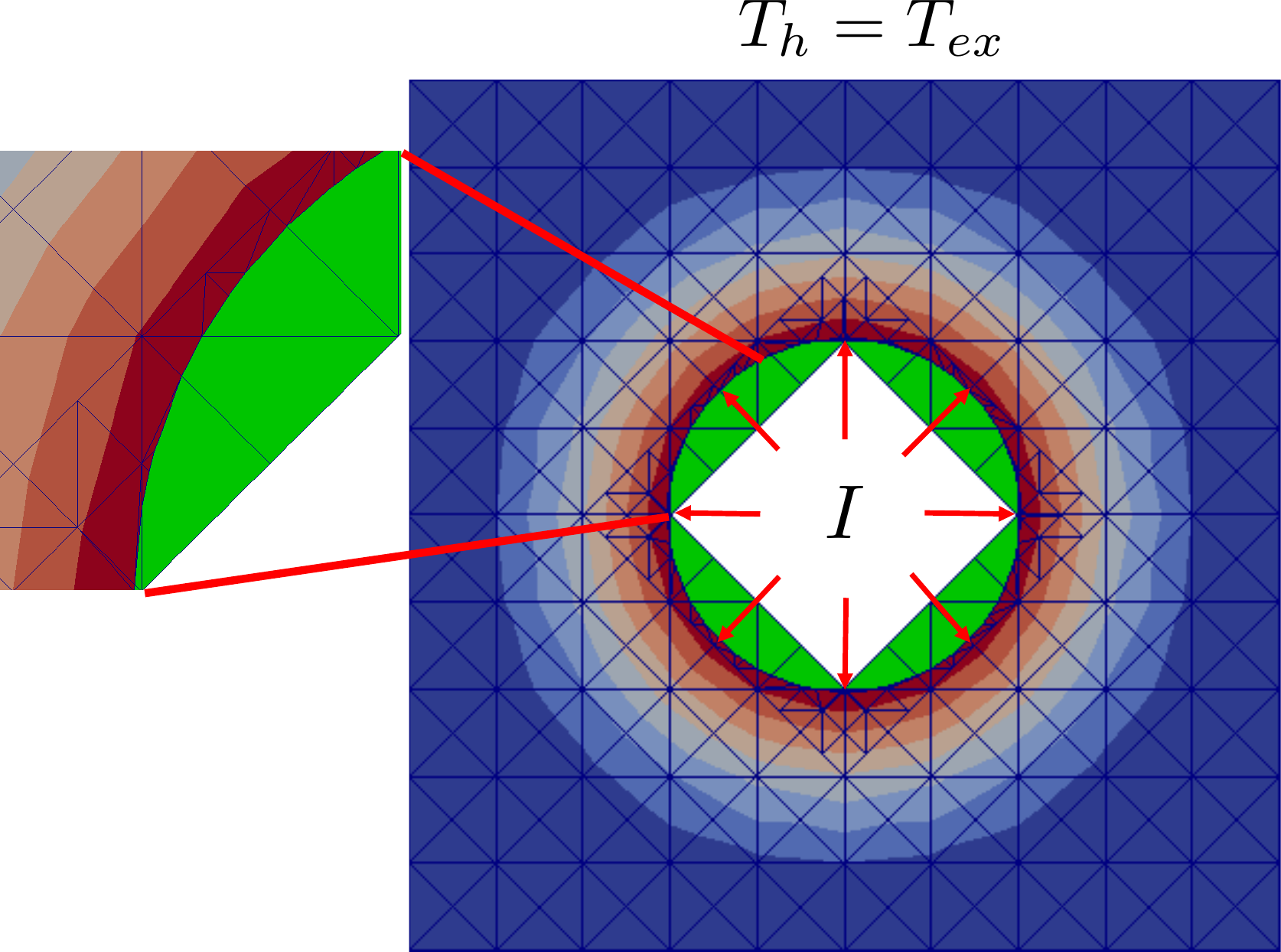}}\vspace{.1cm}
\subfloat[$t=0$.]{\includegraphics[width=.4\textwidth]{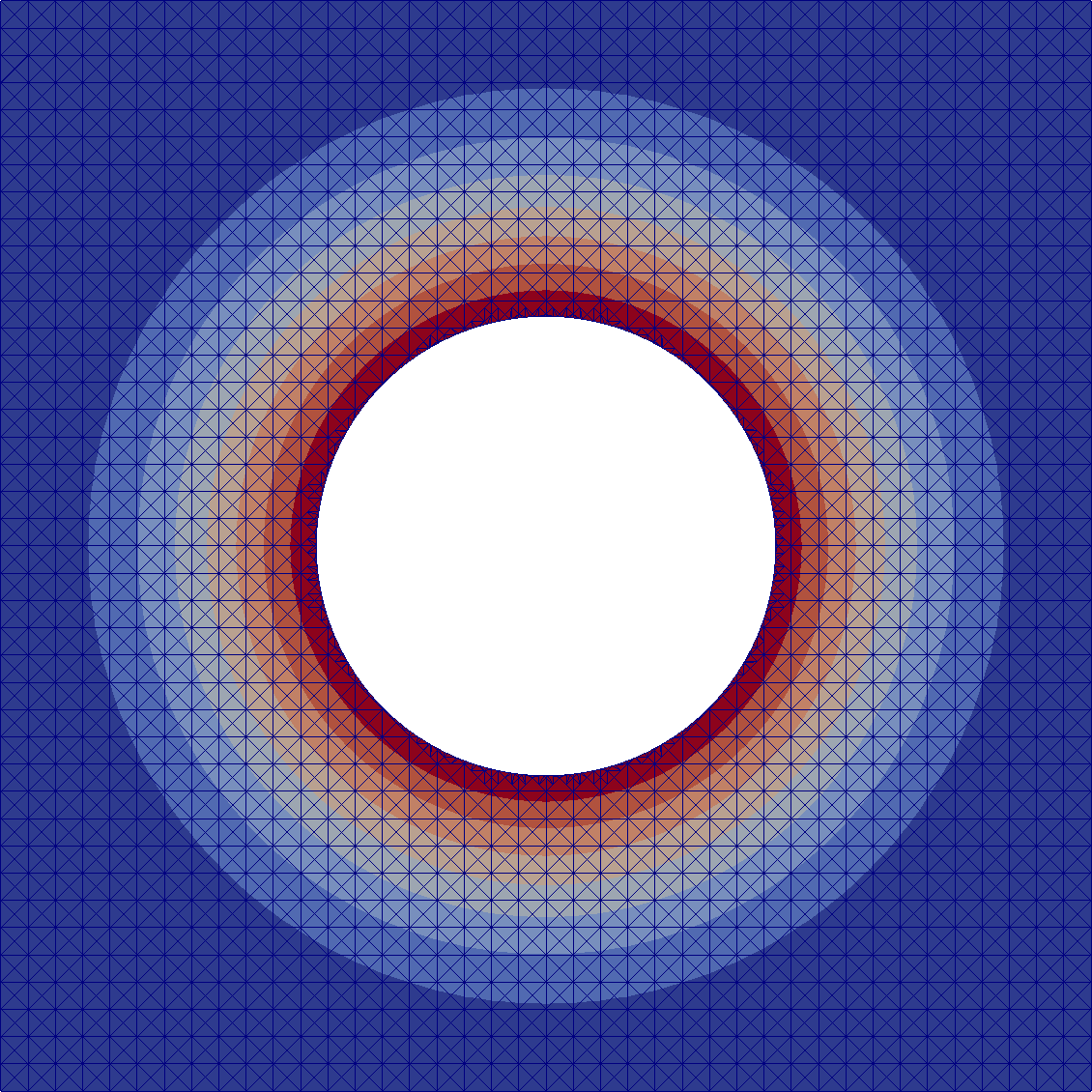}}\\
\subfloat[$t=0.1$.]{\includegraphics[width=.4\textwidth]{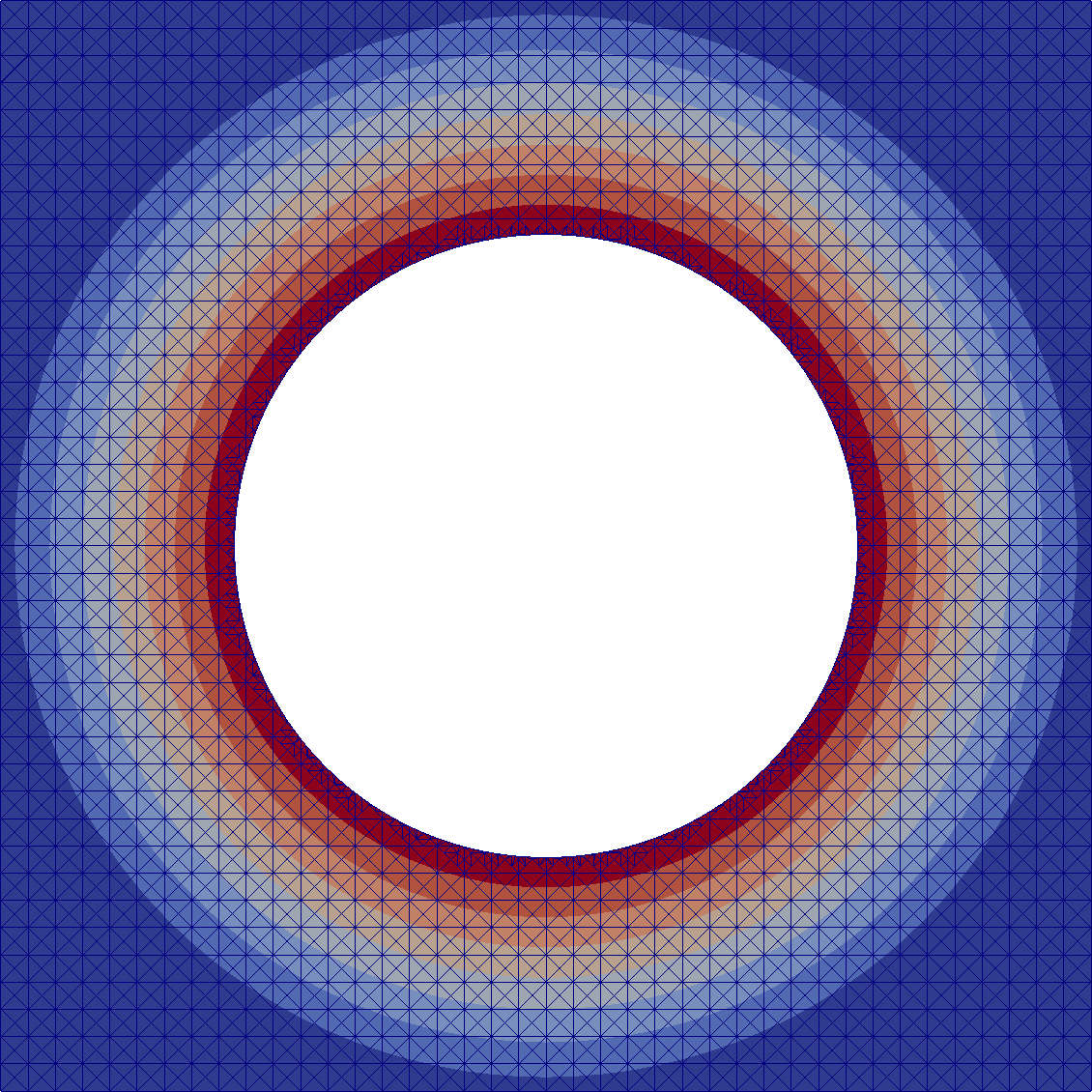}}\vspace{.1cm}
\subfloat[$t=0.2$.]{\includegraphics[width=.4\textwidth]{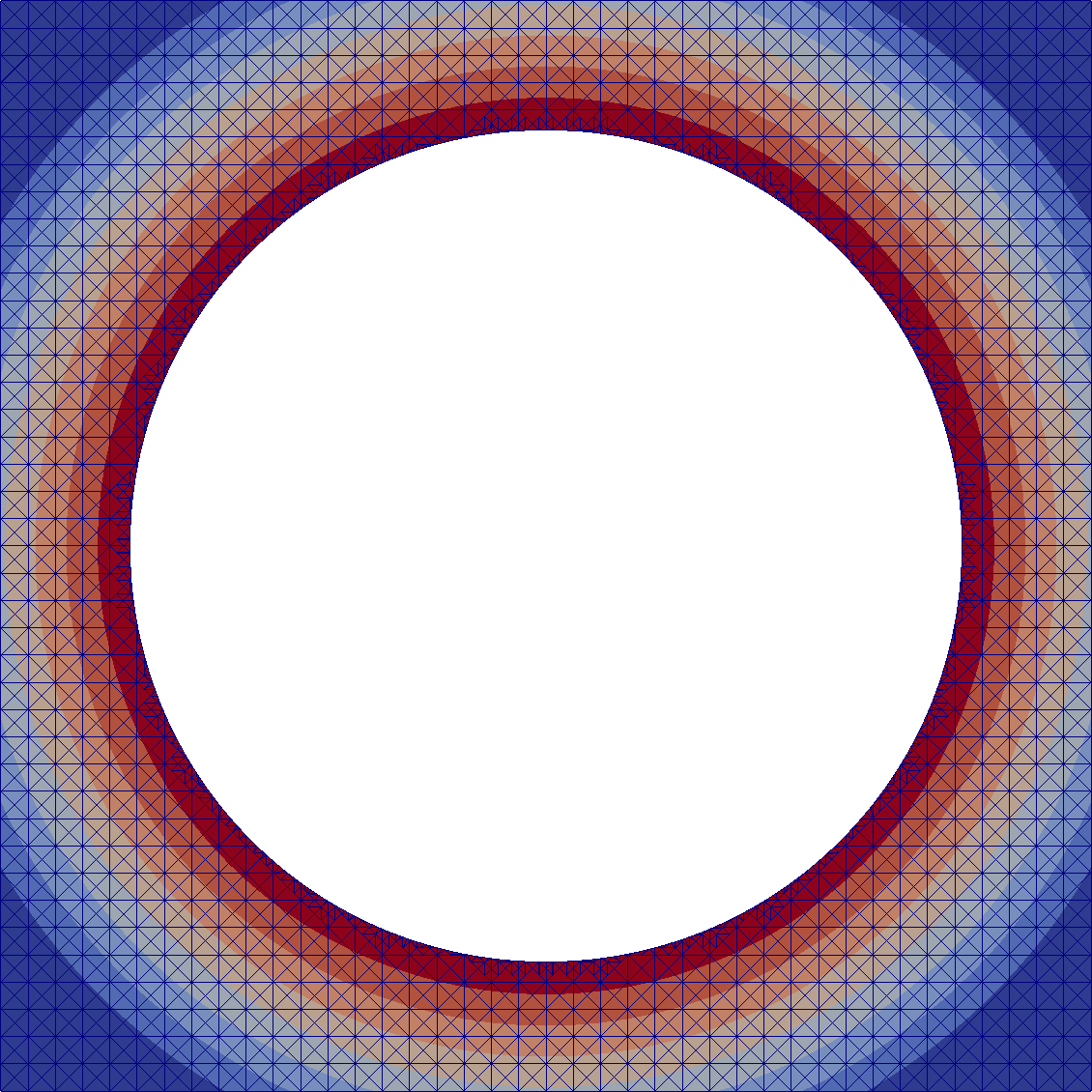}}\vspace{.1cm}
\subfloat{\includegraphics[width=.1\textwidth]{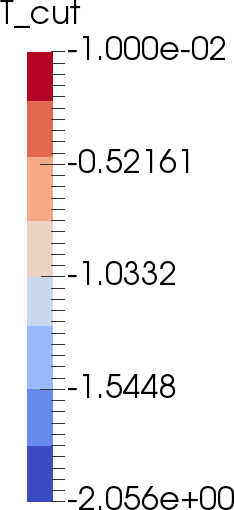}}
\end{center}
\caption{Schematics of the manufactured solution and numerical solution at time $t=\{0,0.1,0.2\}$ for $h=1/40$.}
\label{fig: temperature manufactured solution}
\end{figure}
\begin{figure}
\subfloat[$L^2$ error  in temperature with time step size for $h=1/160$ at time $t=0.1$.]{\includegraphics[width=.5\textwidth]{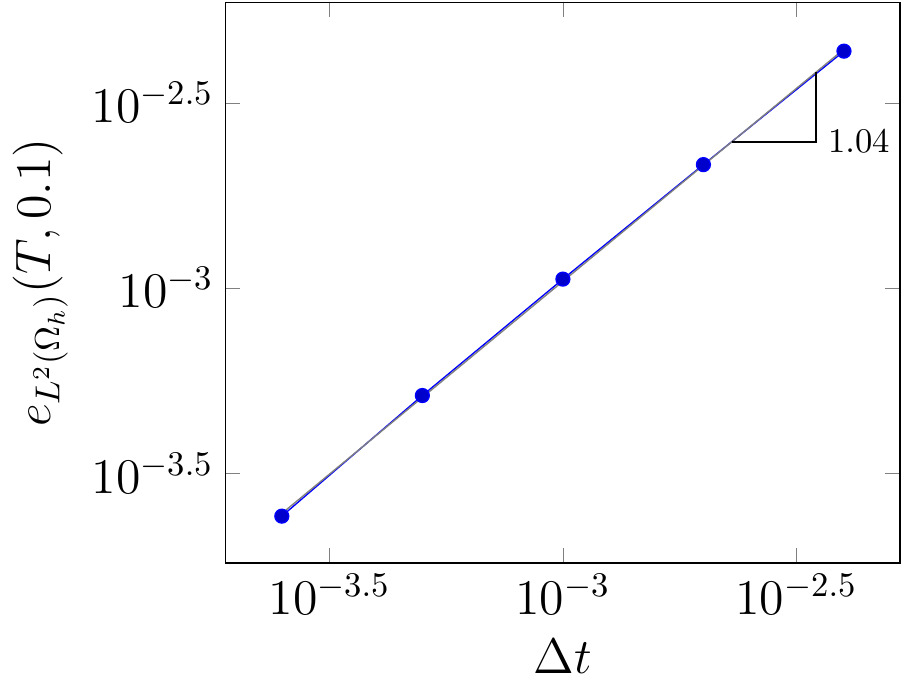}}
\subfloat[$L^2$ error in velocity and radius.]{\includegraphics[width=.5\textwidth]{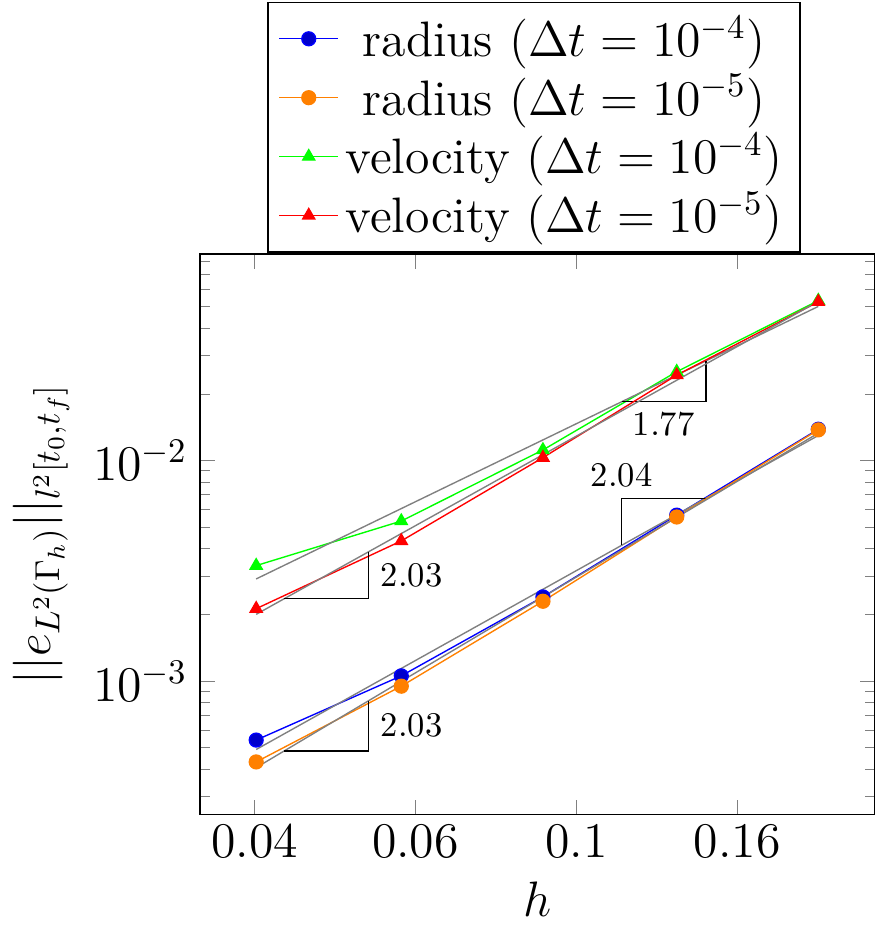}}\\
\subfloat[Temperature $L^2$-error.]{\includegraphics[width=.5\textwidth]{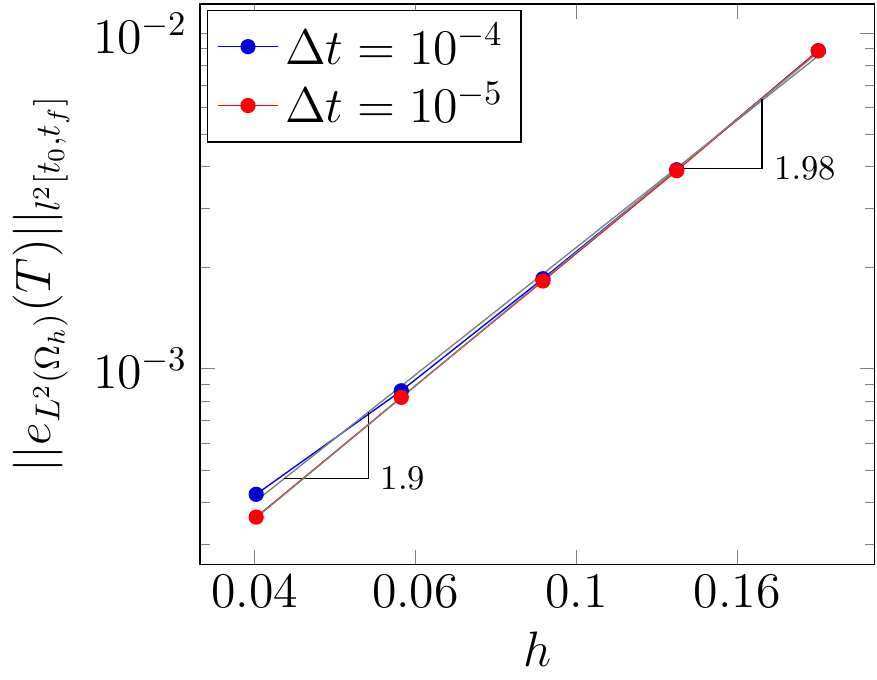}}
\subfloat[Temperature $H^1$-error.]{\includegraphics[width=.5\textwidth]{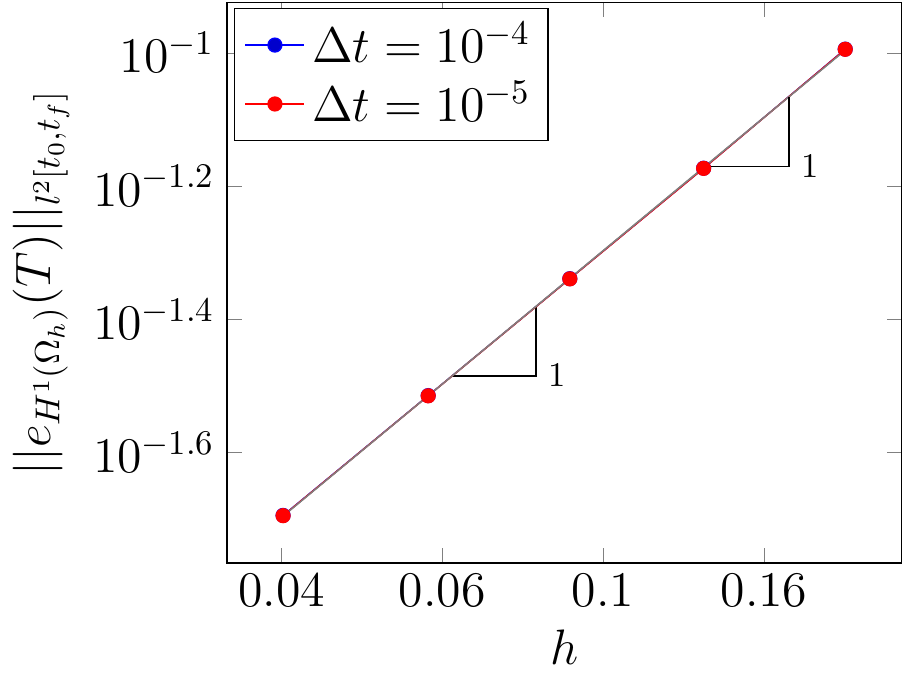}}
\caption{Convergence rates for $L^2$ and $H^1$ errors with mesh refinement and with time step refinement.}
\label{fig: convergence rates for manufactured solution}
\end{figure}
\begin{figure}
\subfloat[Velocity.]{\includegraphics[width=.5\textwidth]{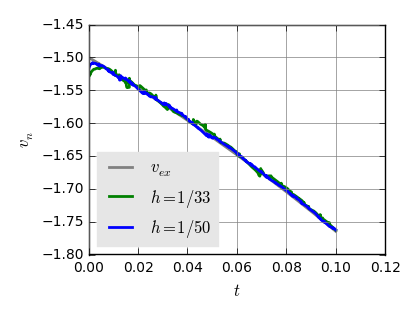}}
\subfloat[Velocity.]{\includegraphics[width=.5\textwidth]{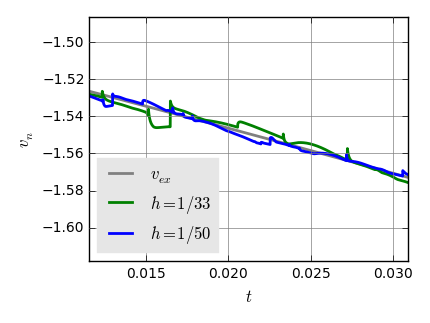}}\\
\subfloat[Radius.]{\includegraphics[width=.5\textwidth]{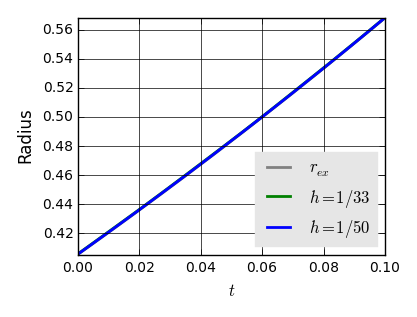}}
\subfloat[Radius.]{\includegraphics[width=.5\textwidth]{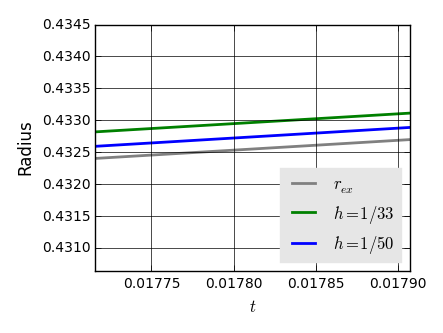}}
\caption{Computed average velocity and average radius versus analytical solution for $\Delta t = 10^{-5}$.}
\label{fig: velocity and radius}
\end{figure}

We have constructed a two-dimensional manufactured solution inspired by the manufactured solution of a  two-phase Stefan problem in \cite{Jahn2017}. We consider a rectangular domain $\Omega$ with a circular hole. The circular hole gets heated by a heat flux $I(x,t)$. We choose $T_m=-0.01$, $\rho=c=k=1.0$, $L=1.0$. We consider an analytical temperature distribution given by  
\begin{equation}
T_{ex}(x,t) = -e^{r(x)} + \cos\left (\frac{\pi r(x) }{2 \log\left(\alpha(t) \right )}\right) - T_m + \alpha(t). 
\end{equation}
Here, $r(x)=\sqrt{\text{x}^2+\text{y}^2}$ and 
\begin{equation}
\alpha(t)=\frac{3}{2-3t}.
\end{equation}
At $r(x)=R(t)$ with 
\begin{equation}
R(t) = \log(\alpha(t))
\end{equation}
the analytical solution is at melting temperature $T=T_m$. For $t=0$, we obtain the initial condition with $\alpha(0)=\frac{3}{2}$ as
\begin{equation}
T_0 = -e^{r(x)} + \cos\left( \frac{\pi r(x)}{2 \log(1.5)}\right) + \frac{3}{2} - T_m.
\end{equation}
The volume source term $f$ can now be determined from $T_{ex}$\footnote{the corresponding symbolic derivation using an IPython notebook can be found in \cite{susanne_claus_2018_1311829}.}. A level-set describing the location of the melting temperature is given by  
\begin{equation}
\phi(x,t) = R(t) - r(x). 
\end{equation}
This level-set describes the motion of the circular hole and the normal velocity of the hole is given by
\begin{equation}
v(x,t) \cdot n_{\Gamma} = - \frac{\partial R(t)}{\partial t} = - \alpha(t).
\end{equation}
The expression for the beam at $\Gamma(t)$ can now be determined from \eqref{equ: latent heat interface} which yields
\begin{equation}
\begin{aligned}
I_{ex}(x,t)  &=  A_{ex}(t) e_{ray,ex}(x), \\
A_{ex}(t ) &= - \left[ (\rho L +1) \alpha(t) + \frac{\pi}{2 R(t)} \right], \\ 
 e_{ray,ex}(x) &= - n_{\Gamma,ex} = \frac{1}{r(x)} \begin{pmatrix} \text{x} \\ \text{y} \end{pmatrix}.
\end{aligned}
\label{equ: beam manufactured solution}
\end{equation}
To test our numerical scheme, we set the temperature $T=T_{ex}$ on $\partial \Omega$ and apply the heat flux expression \eqref{equ: beam manufactured solution} on $\Gamma(t)$. Figure~\ref{fig: temperature manufactured solution} shows the numerical solution at time $t={0,0.1,0.2}$ and shows the removal of material with time. We test convergence with mesh refinement and with time step refinement.  We evaluate the error of our numerical solution with respect to the analytical solution in the following relative error norms. For each time-step, $t_n \in {t_0, ..., t_{n_t}}$, we determine 
\begin{equation}
\begin{aligned}
e_{L^2(\Omega_h)}(u,t_n) := \frac{||u_h - u_{ex} ||_{L^2(\Omega_h(t_n))}}{||u_{ex} ||_{L^2(\Omega_h(t_n))}} = \frac{\sqrt{ \int_{\Omega_h(t_n)} (u_h-u_{ex})^2 \, dx}}{\sqrt{ \int_{\Omega_h(t_n)} (u_{ex})^2 \, dx}}, \\
e_{H^1(\Omega_h)}(u,t_n) := \frac{||u_h - u_{ex} ||_{H^1(\Omega_h(t_n))}}{||u_{ex} ||_{H^1(\Omega_h(t_n))}} = \frac{\sqrt{ \int_{\Omega_h(t_n)} (u_h-u_{ex})^2 + (\nabla(u_h-u_{ex}))^2 \, dx}}{\sqrt{ \int_{\Omega_h(t_n)} (u_{ex})^2 + \nabla u_{ex}^2 \, dx}}, \\
e_{L^2(\Gamma_h)}(u,t_n) :=\frac{||u_h - u_{ex} ||_{L^2(\Gamma_h(t_n))}}{||u_{ex} ||_{L^2(\Gamma_h(t_n))}} = \frac{\sqrt{ \int_{\Gamma_h(t_n)} (u_h-u_{ex})^2 \, dx}}{\sqrt{ \int_{\Gamma_h(t_n)} (u_{ex})^2 \, dx}}. \\
\end{aligned}
\label{equ: error norms}
\end{equation}
Here, $u_h$ is the numerical solution and $u_{ex}$ is the analytical solution, which in the following will be the temperature, the radius or the interface velocity. To average the error over time, we define the $l^2$-error norm over the time interval as 
\begin{equation}
\begin{aligned}
||e(u)||_{l^2[t_0,t_f]} = \sqrt{ \frac{1}{n_t} \sum_{i=0}^{n_t} e(u,t_i)^2}.
\end{aligned}
\end{equation}
Here, $e(u)$ is any of the error measures defined in equation~\eqref{equ: error norms}. \\
We compute the numerical solution in the time interval $t \in [0,0.1]$ for time step sizes $\Delta t = \{ 10^{-4},10^{-5}\}$. Figure~\ref{fig: convergence rates for manufactured solution} shows the convergence of the temperature to the analytical solution with mesh refinement and the convergence of the temperature with time step refinement. As can be seen in Figure~\ref{fig: convergence rates for manufactured solution}, we obtain optimal convergence orders of second order for the $L^2$-norm and of first order for the $H^1$-norm in space and first order convergence in the $L^2$-norm in time for the temperature. For the convergence of velocity and radius of the circular hole, we obtain a convergence rate of second order. The convergence rates for the $L^2$-errors in temperature, velocity and radius show a slight improvement for the time step $\Delta t = 10^{-5}$ in comparison to time step $\Delta t = 10^{-4}$ for finer meshes. This is to be expected as the discretisation error in time is starting to dominate the total error for finer meshes and is impacting the convergence rate. Figure~\ref{fig: velocity and radius} shows the averaged computed velocity and the averaged computed radius and their analytical expression. The average is computed for all $t_n \in {t_0, ..., t_{n_t}}$ as
\begin{equation}
\begin{aligned}
v_{avg}(t_n) &= \frac{\int_{\Gamma_h(t_n)} v_n^h \, ds}{\int_{\Gamma_h(t_n)} ds}, \\
r_{avg}(t_n) &= \frac{\int_{\Gamma_h(t_n)} \sqrt{x^2+y^2} \, ds}{\int_{\Gamma_h(t_n)} ds}.
\end{aligned}
\end{equation}
It is clear that both the velocity and radius approach the exact solution as the mesh is refined. The convergence in position seems to be monotonic, at any time of the analysis, while the instantaneous convergence in velocity appears to be much more erratic, which is to be expected, given the fact that the velocity is the time derivative of the position.

\subsection{Thermal ablation using a moving laser beam}
In this Section, we will set-up several numerical examples describing a laser beam heating a workpiece alongside a predefined machining path. We define the intensity of the spatially Gaussian-distributed beam as 
\begin{equation}
\begin{aligned}
I(x,t) &= - A_p(\theta) f(x,t)  e_{ray}, \\
f(x,t) &=  
f_x(x,t) \, f_t(t) 
:=
A_{amp}\frac{1}{\sqrt{2\pi^{d-1}\sigma^2} }e^{\frac{-p(x,t) \cdot p(x,t)}{2\sigma^2}} \, f_t(t),   \\ 
 p(x,t) &= (x - F(t)) - ((x-F(t)) \cdot e_{ray}) \, e_{ray},
\end{aligned}
\label{equ: Gaussian beam}
\end{equation}
where $\sigma$ is the width of the beam,  $A_{amp}$ is the amplitude of the beam, $F(t)$ is the focal point of the beam that describes the path of the laser beam. In the following, we choose the direction of the beam, $e_{ray}$, to be constant in time. The beam is scaled with the absorption coefficient (\cite{Schulz1987}, \cite{Stratton1941}, \cite{Schulz2017}) given by 
\begin{equation}
A_p(\theta) = \begin{cases}
1 - \frac{2\cos(\theta)^2 - 2 \varepsilon \cos(\theta) + \varepsilon^2}{2\cos(\theta)^2 + 2\varepsilon\cos(\theta)+\varepsilon^2} \quad& \cos(\theta) > 0, \\
0 \quad &\mbox{otherwise} \, ,
\end{cases}
\end{equation}
where the angle of incidence $\theta$ of the laser beam with respect to the inside surface normal $-n_{\Gamma}$ appears in the equation through trigonometric function $\cos(\theta) = - n_{\Gamma}(x,t) \cdot e_{ray}$.  Here, $\varepsilon$ is a material-dependent quantity, which we choose as $\varepsilon=1$. We choose to represent a pulsed laser beam whose periodic on/off behaviour can described by using the pulse function 
\begin{equation}
f_t(t) = 
\begin{cases}
1 &  \mbox{if } t - \left\lfloor{\frac{t}{P_0} }\right\rfloor P_0  \leq \frac{P0}{2}, \\
0  &\mbox{else},
\end{cases}
\end{equation}
where $ \lfloor \ \rfloor$ denotes the floor operation, and $P_0$ is the total period, which is the sum of an "on" phase of duration  $t_{ON}$ and an "off" phase of duration $t_{OFF}$ during which the workpiece does not receive any energy from the thermal ablation device. \\
In the following sections, we consider rectangular workpieces for which the top boundary is the moving boundary $\Gamma$ (i.e. thermally ablated surface).  Homogenous Dirichlet boundary conditions, i.e. $T |_{\partial \Omega_D} =0$, will be applied to the bottom boundary, and homogeneous Neumann boundary conditions will be applied to the remaining sides (see Figure~\ref{fig: laser cutting schematics}).

\subsubsection{Pulsed thermal ablation in 2D}

\begin{figure}
\begin{center}
\subfloat[$t=0.4, P_0=0.01$.]{\includegraphics[width=.45\textwidth]{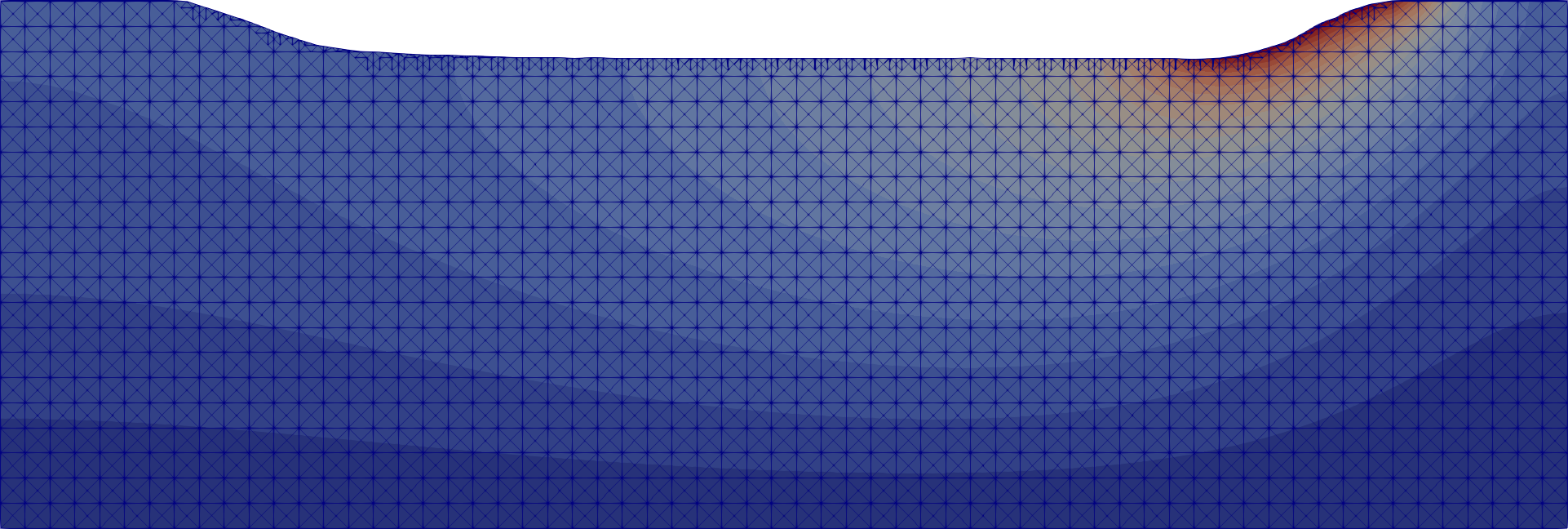}} \vspace{.1cm}
\subfloat[$t=0.4, P_0=0.1$.]{\includegraphics[width=.45\textwidth]{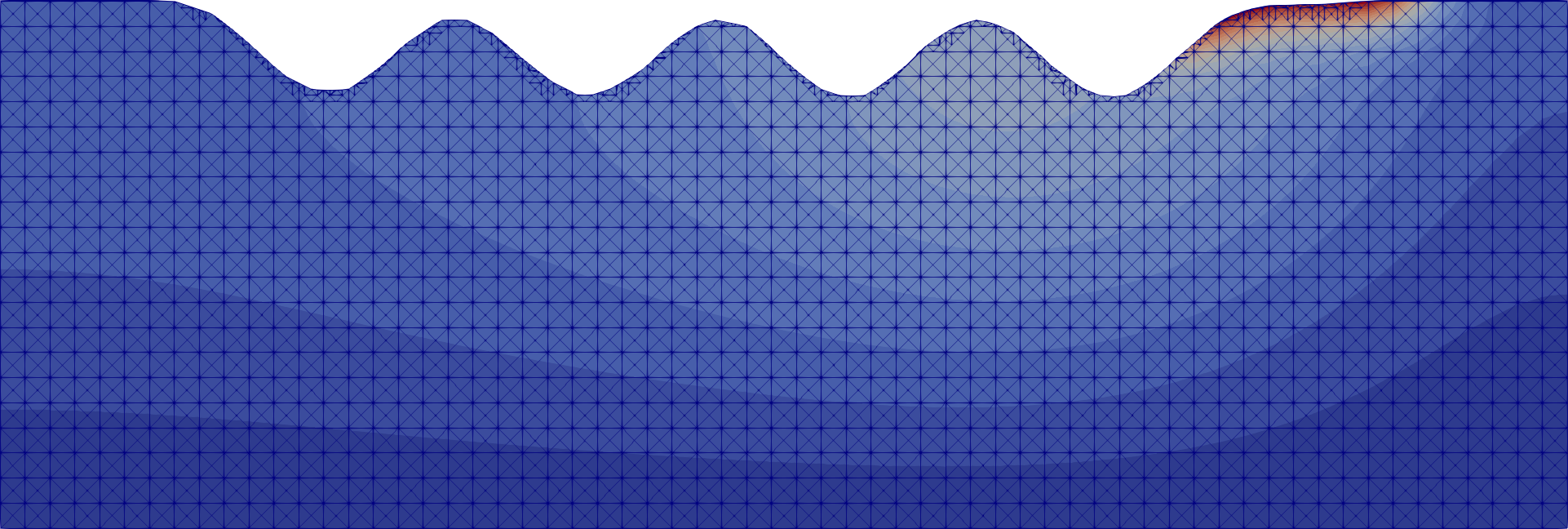}}\\ 
\subfloat[$t=0.8, P_0=0.01$.]{\includegraphics[width=.45\textwidth]{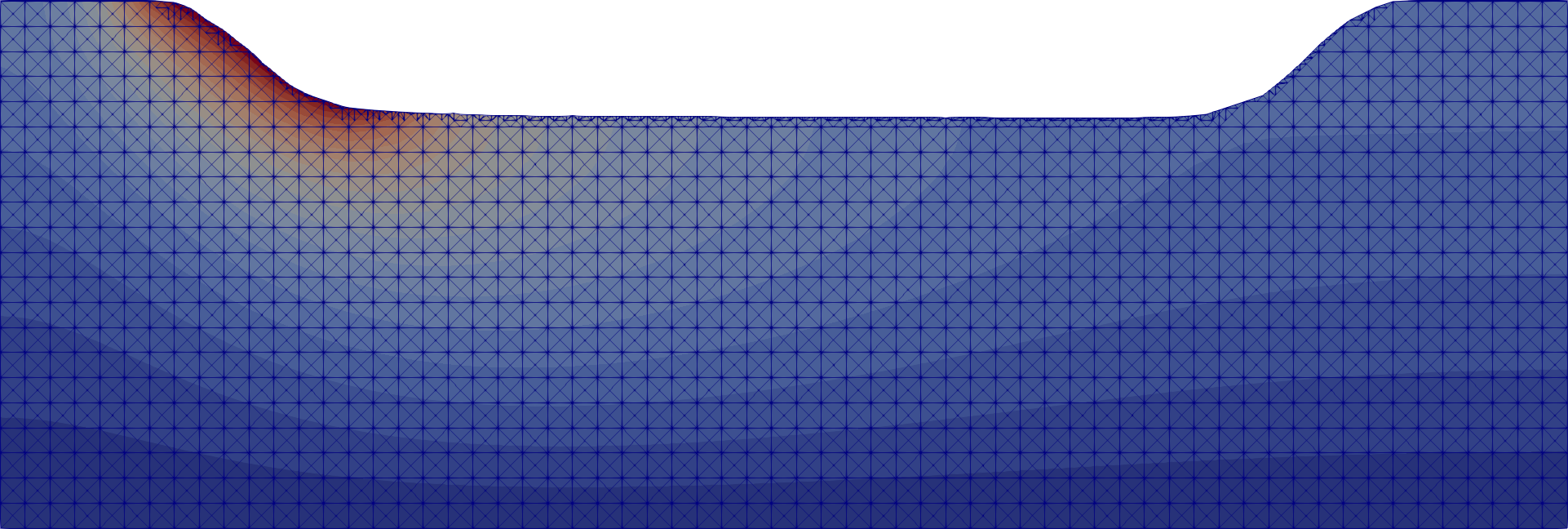}}\vspace{.1cm}
\subfloat[$t=0.8, P_0=0.1$.]{\includegraphics[width=.45\textwidth]{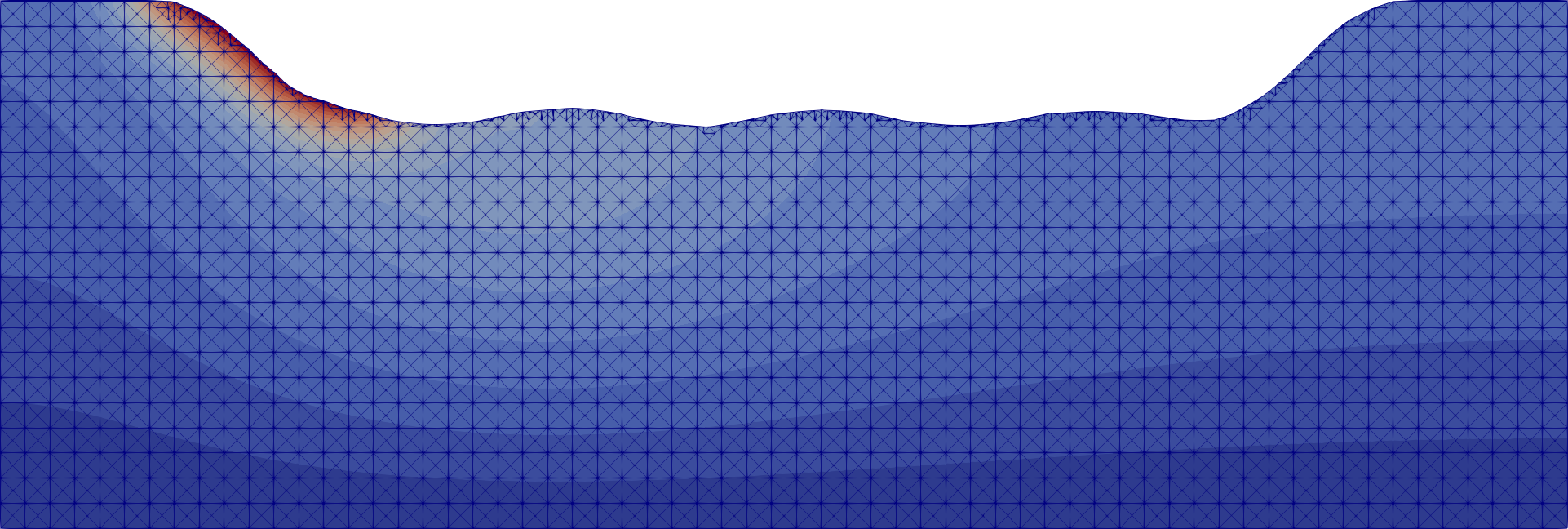}}\\
\subfloat[$t=1.2, P_0=0.01$.]{\includegraphics[width=.45\textwidth]{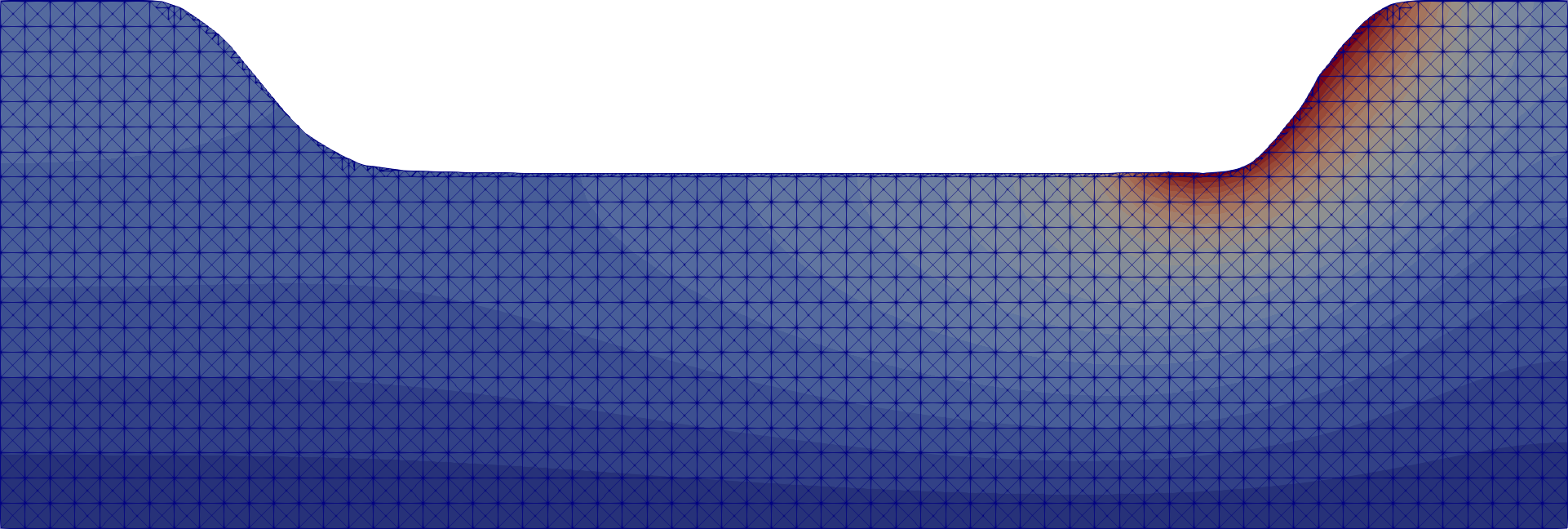}}\vspace{.1cm}
\subfloat[$t=1.2, P_0=0.1$.]{\includegraphics[width=.45\textwidth]{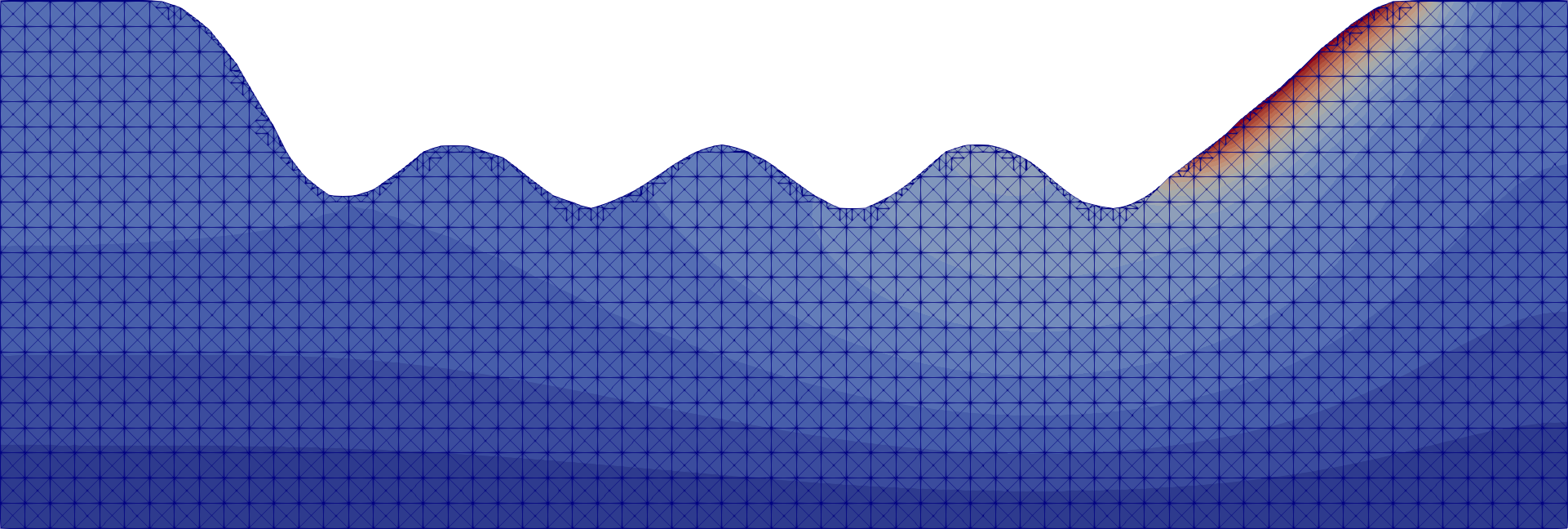}}\\
\subfloat[$t=1.6, P_0=0.01$.]{\includegraphics[width=.45\textwidth]{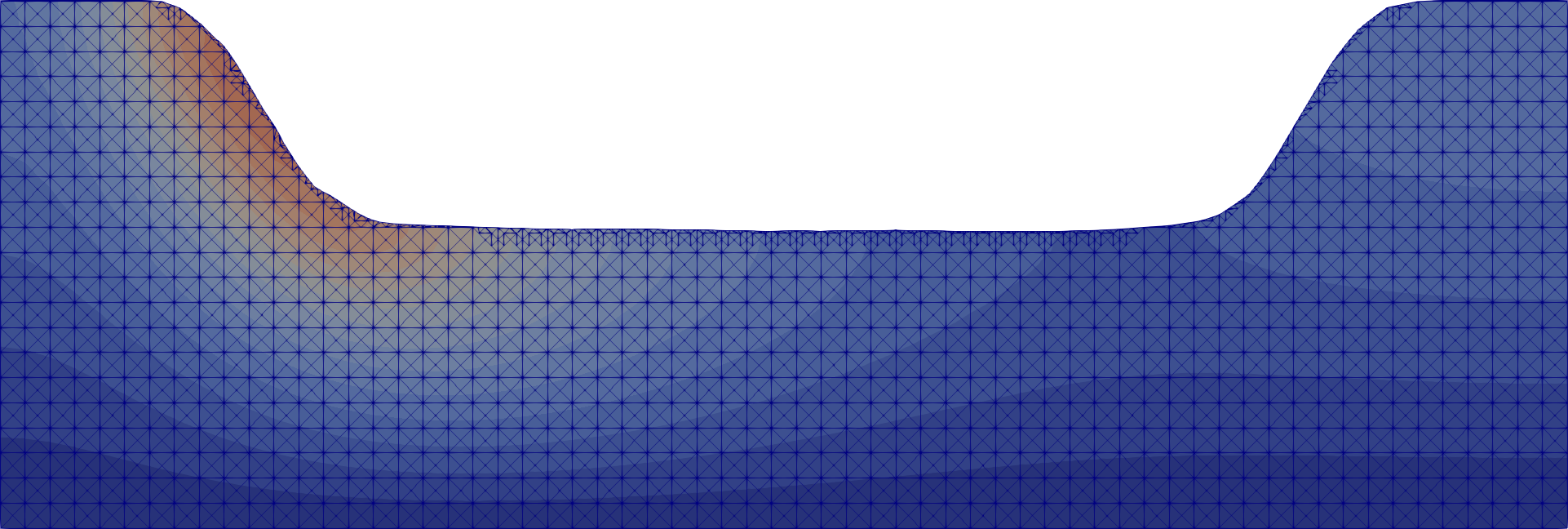}}\vspace{.1cm}
\subfloat[$t=1.6, P_0=0.1$.]{\includegraphics[width=.45\textwidth]{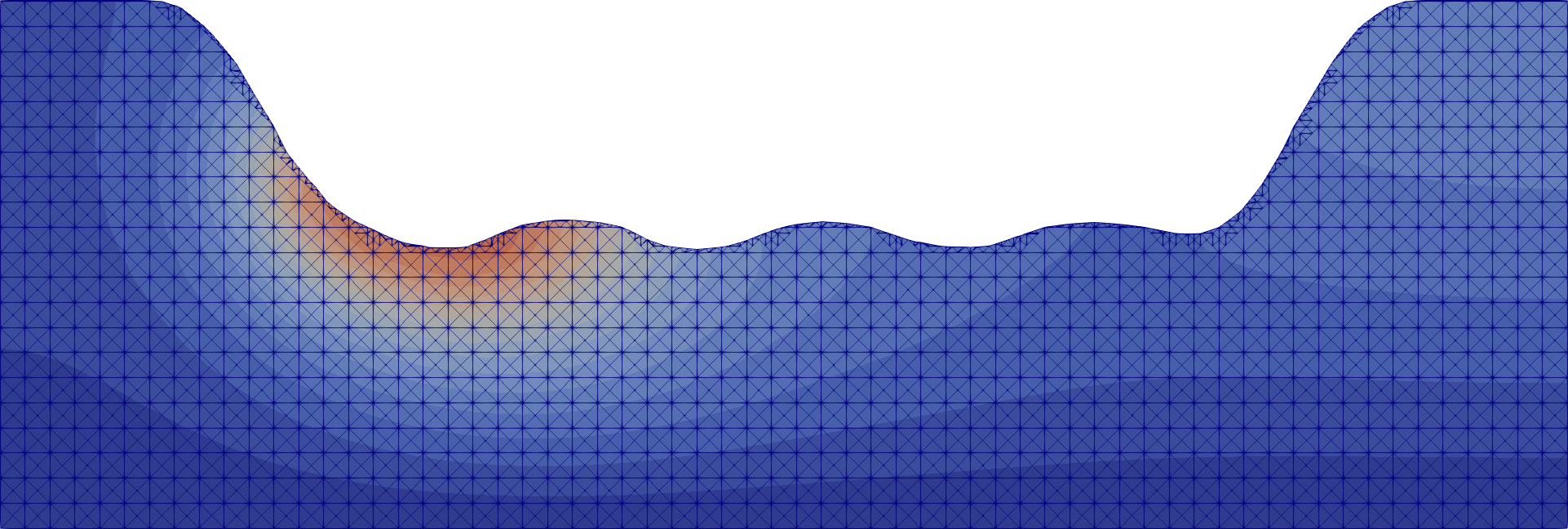}}\\
\subfloat{\includegraphics[width=.45\textwidth]{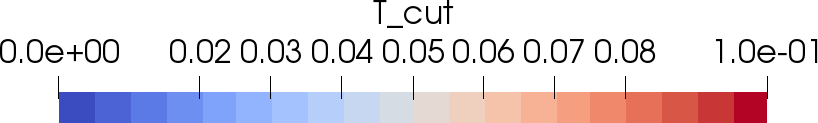}}
\end{center}
\caption{Pulsed laser beam for periods $P_0=0.01$ on the left and $P_0=0.1$ on the right for time $t=0.4,0.8,1.2,1.6$.}
\label{fig: temperature laser 2D}
\end{figure}

\begin{figure}
\begin{center}
\includegraphics[width=.7\textwidth]{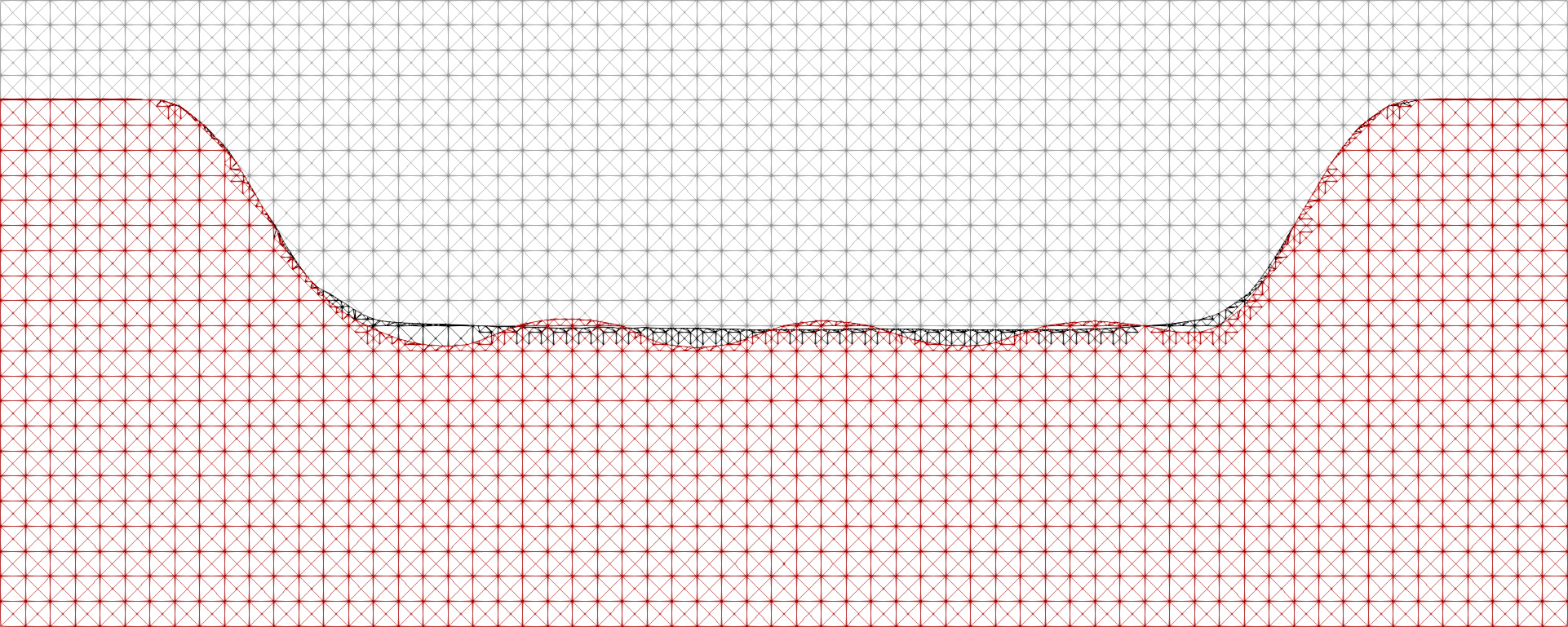}
\end{center}
\caption{Background mesh and cut-meshes at final time $t=1.6$ for pulsed beam $P_0=\{0.01,0.1\}$.}
\label{fig: temperature laser 2D cutmesh}
\end{figure}

Consider a rectangular background domain $\Omega_b=(0,3)\times(0,1.2)$ and a time interval $t \in [0,1.6]$. We consider an initial level-set of $\phi(x,0) = y-1.0$ leaving a rectangular block of material $\Omega_h = (0,3)\times(0,1)$. The workpiece $\Omega_h$ is heated by a laser beam described by equation~\eqref{equ: Gaussian beam} with width $\sigma=0.1$, amplitude $A_{amp}=2$ and beam direction $e_{ray} = ( 0 , -1 )$. The time evolution of the laser beam is described by the path of focal point $F(t)$, which is defined for all times $\{t_n\}_{n \in \jump{0\,n_t-1}}$ by
\begin{equation}
 F(t_{n+1}) = F(t_n) + v_F(t) \Delta t
\end{equation}
with initial beam focal point $F_0 = (0.5 \ 1)^T$ and initial velocity $v_F(t)=(5 \ 0)^T$. Velocity $v_F$ conserves its magnitude throughout the simulation, but changes direction every $t_{change} = 0.4$ units of time.  Changing the sign of the velocity vector causes the laser beam to pass over the block of material four times. We choose the material parameters as $T_m = 0.1$, $L=1$, $\rho=1$, $k=1$, $c=1$. We choose two different pulse periods $P_0 \in \{ 0.1 , 0.01 \}$ and compare the corresponding results.  We choose a fixed time step of $\Delta t = 5 \cdot 10^{-4}$ and a fixed mesh size $h=0.048$. Figure~\ref{fig: temperature laser 2D} shows the temperature contour at times $t=0.4,0.8,1.2,1.6$. The short pulsed beam ($P_0=0.01$) removes material in an even manner, leaving no visible crater on the surface of the workpiece, while the long pulsed beam ($P_0=0.1$) leaves a wavy surface with visible craters.
For these two simulations, the amplitude and spatial distribution of the energy is the same, and the ratio between the on and off time are also equal. As a result, the average power received by the workpiece over one period is the same in both cases, which explains why the depths of the resulting cavities are similar (see Figure~\ref{fig: temperature laser 2D cutmesh}). Note however that there is no theoretical reason for a strict equality between the volumes removed during the process as the amount of energy lost through the Dirichlet conditions and the quantity of thermal energy remaining in the workpiece at the end of the simulation may differ (slightly) in the two examples. 

\subsubsection{Laser beam in 3D}
In this Section, we consider two 3D examples. The first example describes the formation of a single crater, for a spatially fixed laser beam, while the second example describes a complex ablation process designed to manufacture a rectangular cavity through the continuous motion of the laser beam. 
 
\begin{figure}
\begin{center}
\subfloat{\includegraphics[width=.7\textwidth]{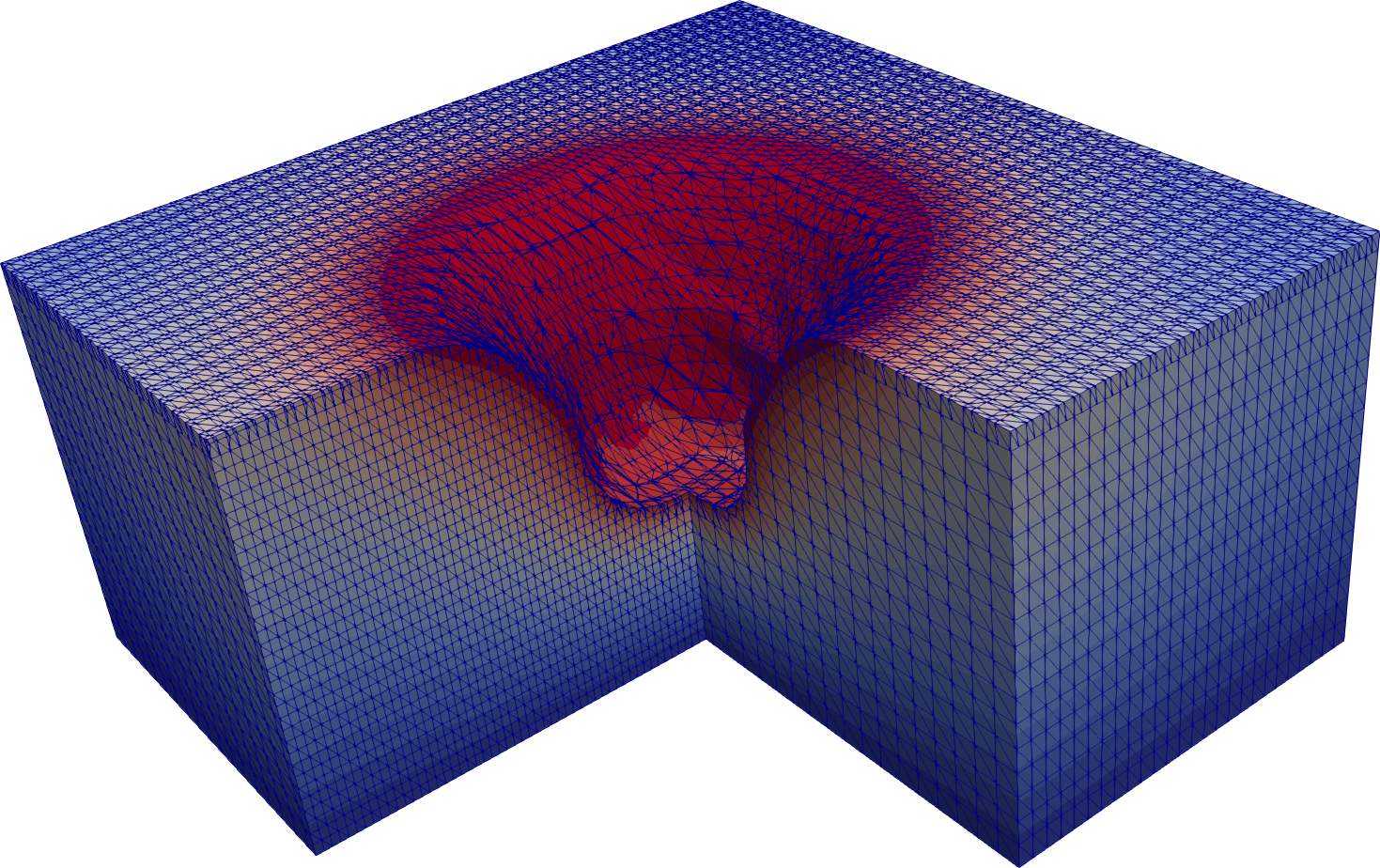}} \vspace{.2cm}
\subfloat{\includegraphics[width=.15\textwidth]{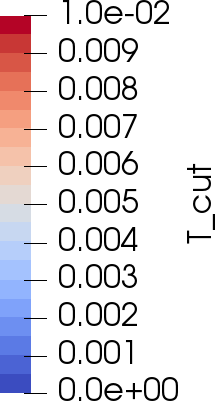}}
\end{center}
\caption{One crater formed by a laser beam at a fixed spatial location, and the distribution of temperature computed at $t=0.2$ units of time.}
\label{fig: one crater}
\end{figure}

\paragraph{Single Crater formation}
We compute the formation of a single crater considering a rectangular background domain $\Omega_b=(-0.5,0.5) \times (-0.5,0.5) \times (-0.5,0.01)$ with an initial level set $\phi(x,0) = y$ in a time interval $t=[0,0.2]$. We fix the time step size to $\Delta t =0.005$ and the mesh size to  $h=0.029$. We choose the material parameters as $T_m = 0.01$, $L=1$, $\rho=1$, $k=1$, $c=1$. The focal point of the laser beam is fixed in time, $F(t) = (0,0,0)$ and we set  $e_{ray} = (0,0,-1)$, $A_{amp}=3$ and 
$\sigma= 0.1$. The laser beam is switched on over the entire time period. Figure~\ref{fig: one crater} shows the crater profile, together with several temperature isolines, at time $t=0.2$. The laser beam causes the formation of a single deep crater in a cone shape. 

\begin{figure}
\begin{center}
\subfloat{\includegraphics[width=.6\textwidth]{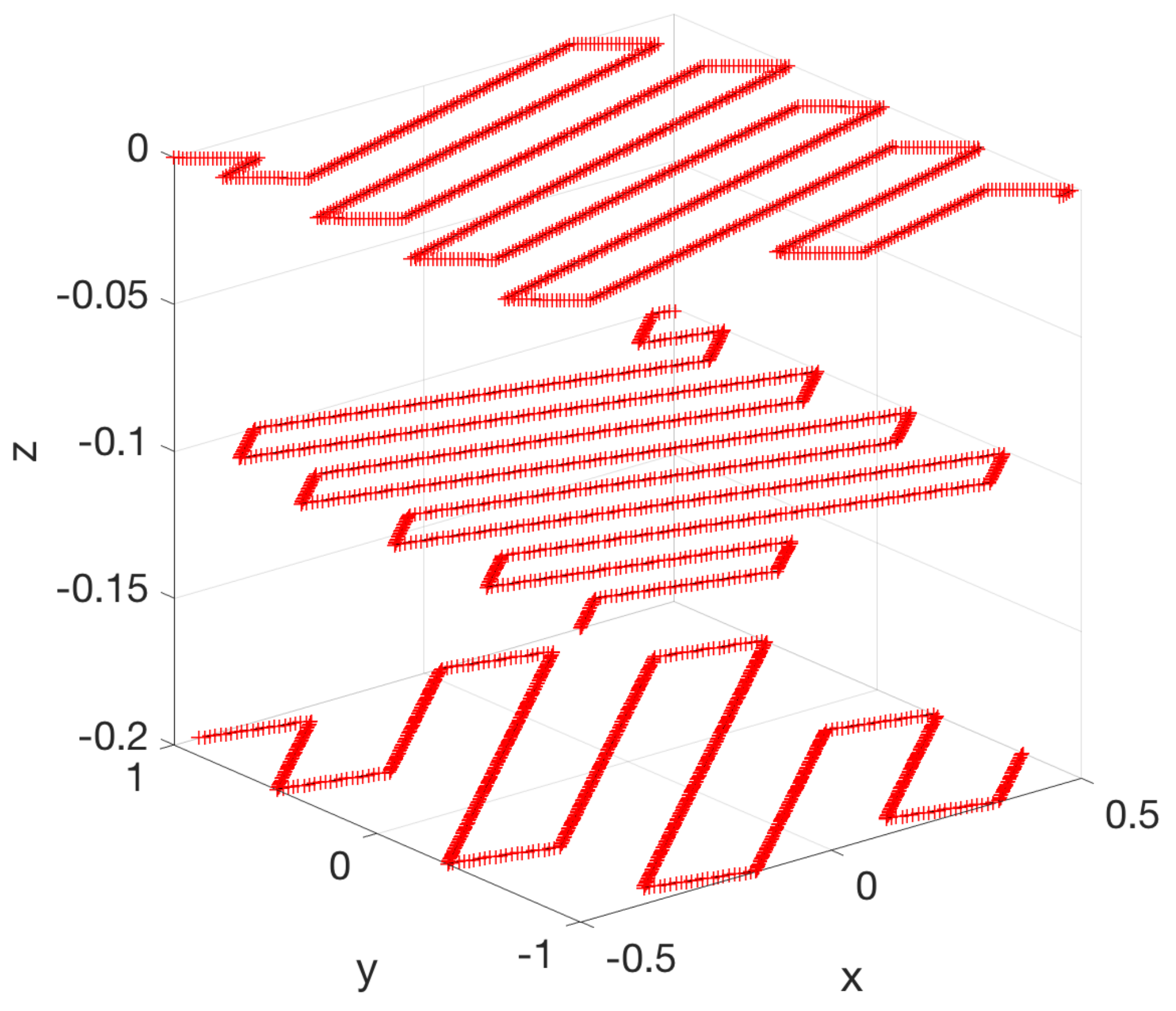}} 
\end{center}
\caption{Prescribed path of the focal point of a laser beam designed to create a rectangular cavity. The depth is not significant here as the laser beam is invariant in the $z$ direction.}
\label{fig: 3D path}
\end{figure}

\begin{figure}
\begin{center}
\subfloat[$t=0$.]{\includegraphics[width=.45\textwidth]{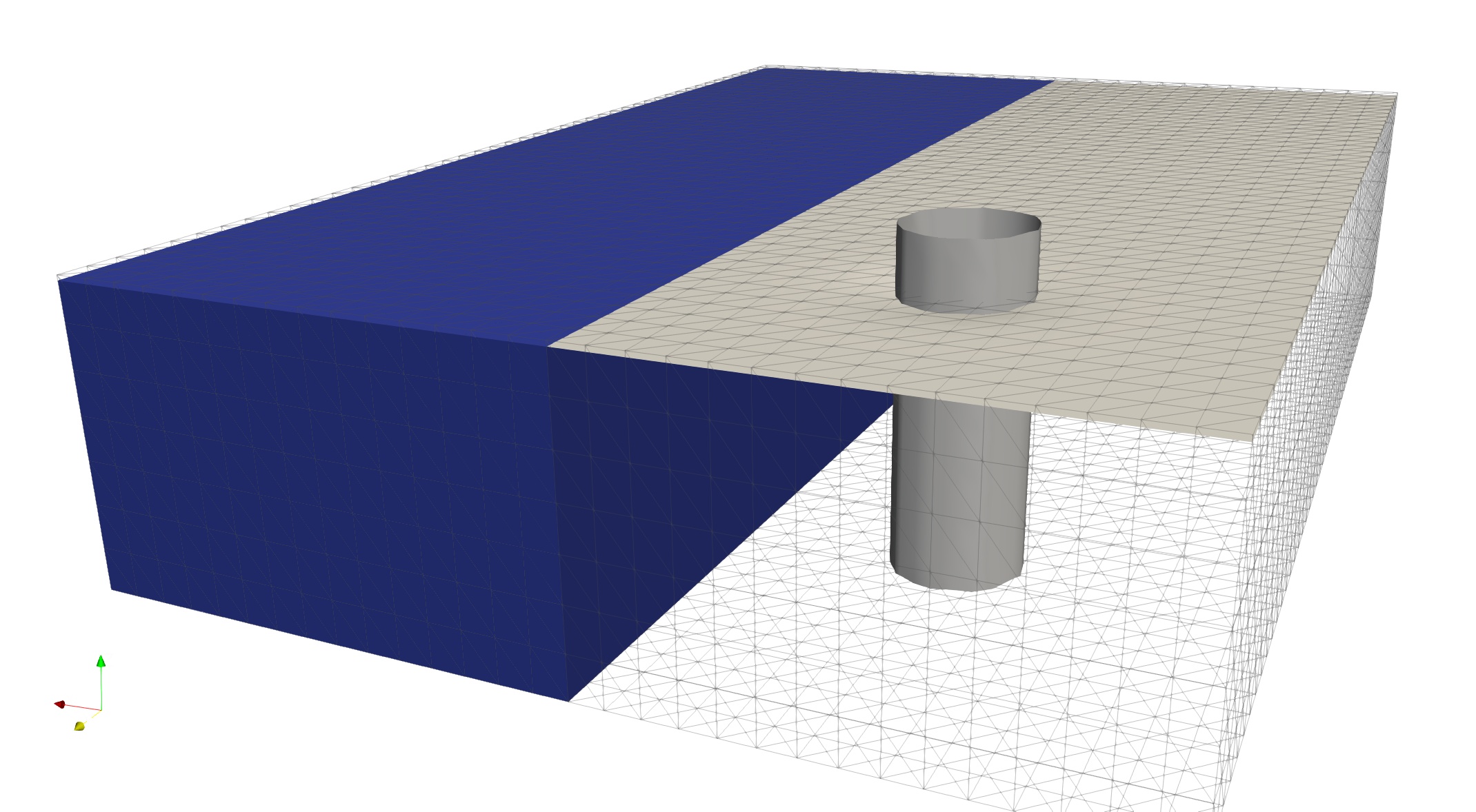}} \subfloat[$t=1$.]{\includegraphics[width=.45\textwidth]{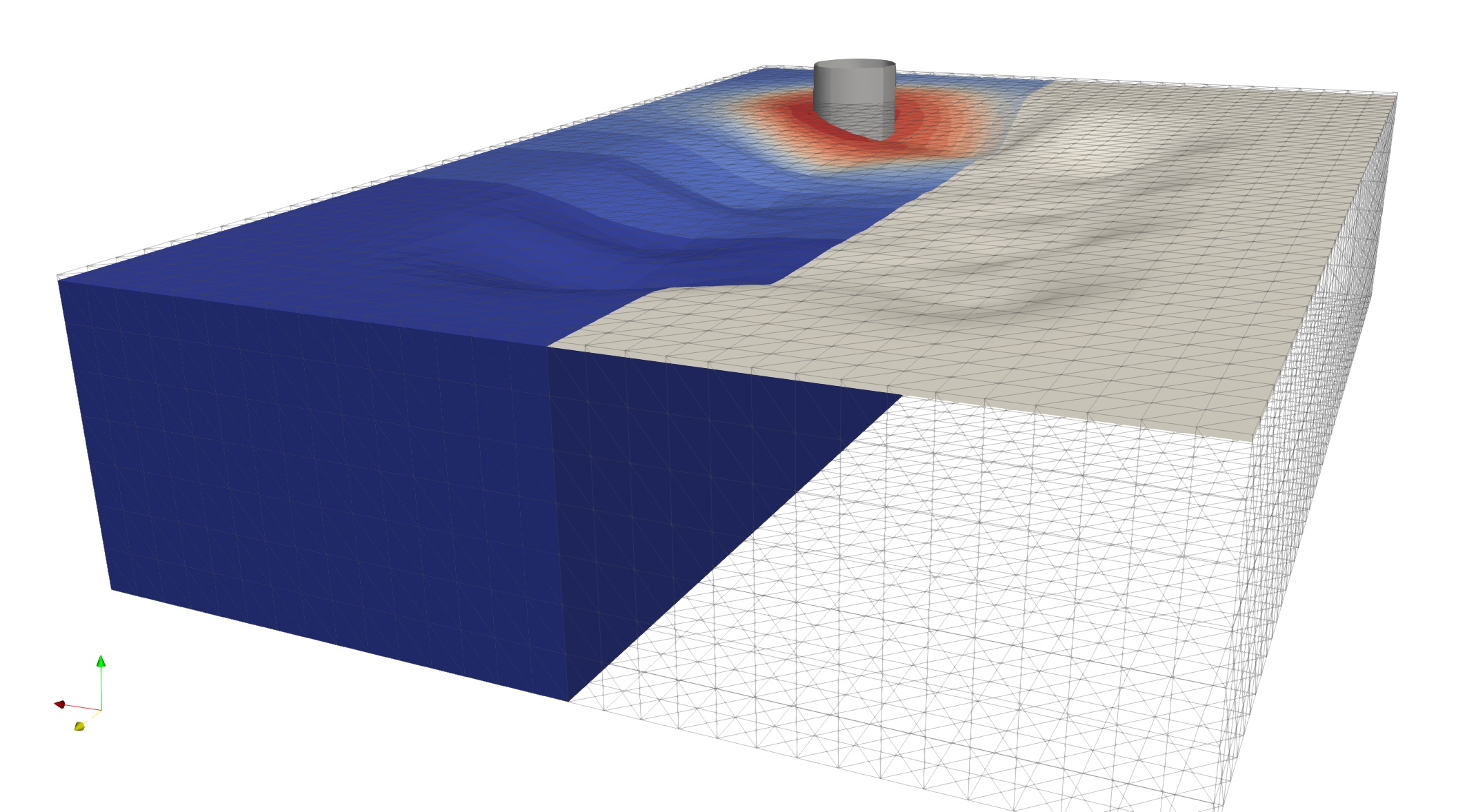}}\vspace{.4cm}
\\ 
\subfloat[$t=2$.]{\includegraphics[width=.45\textwidth]{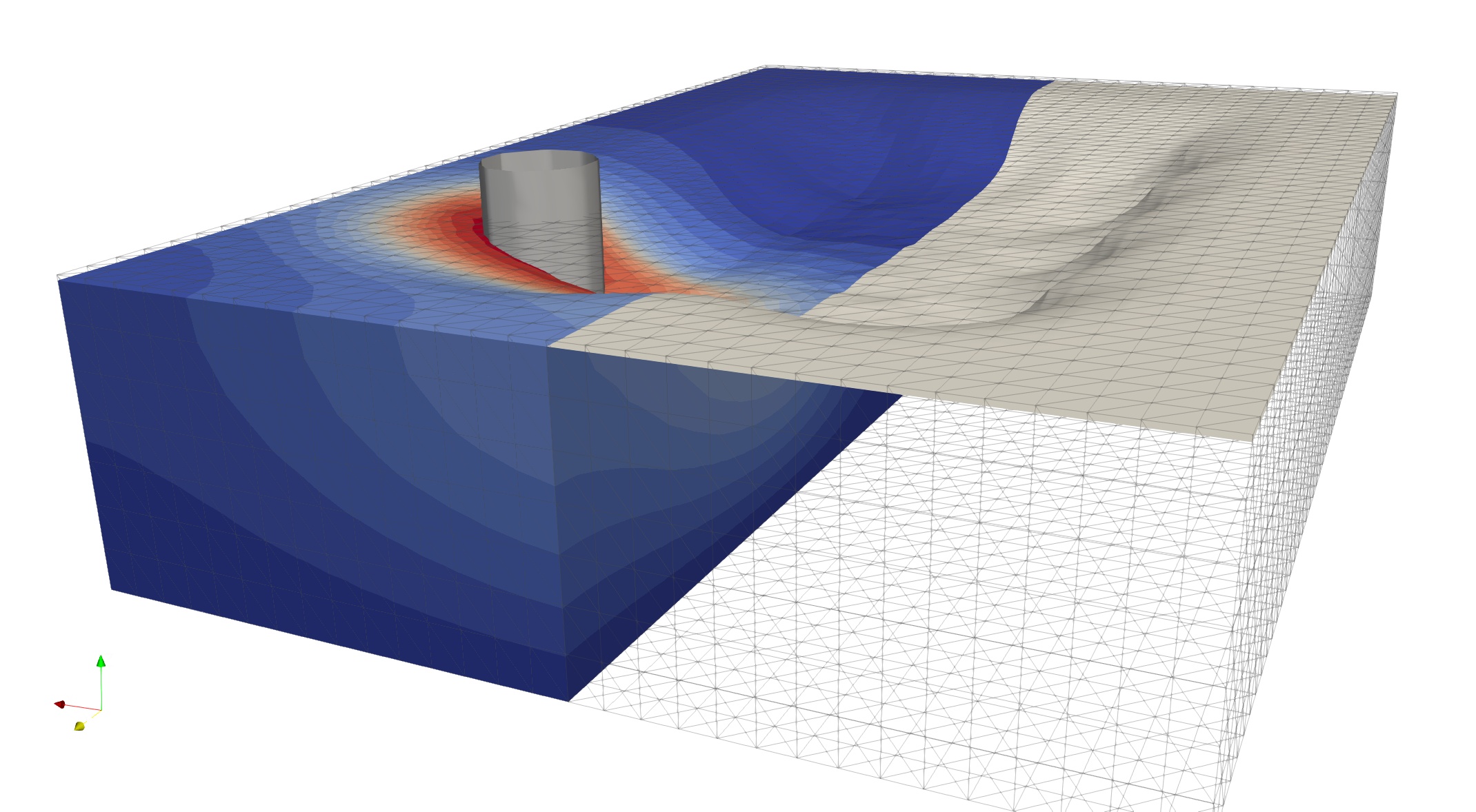}}
\subfloat[$t=3$.]{\includegraphics[width=.45\textwidth]{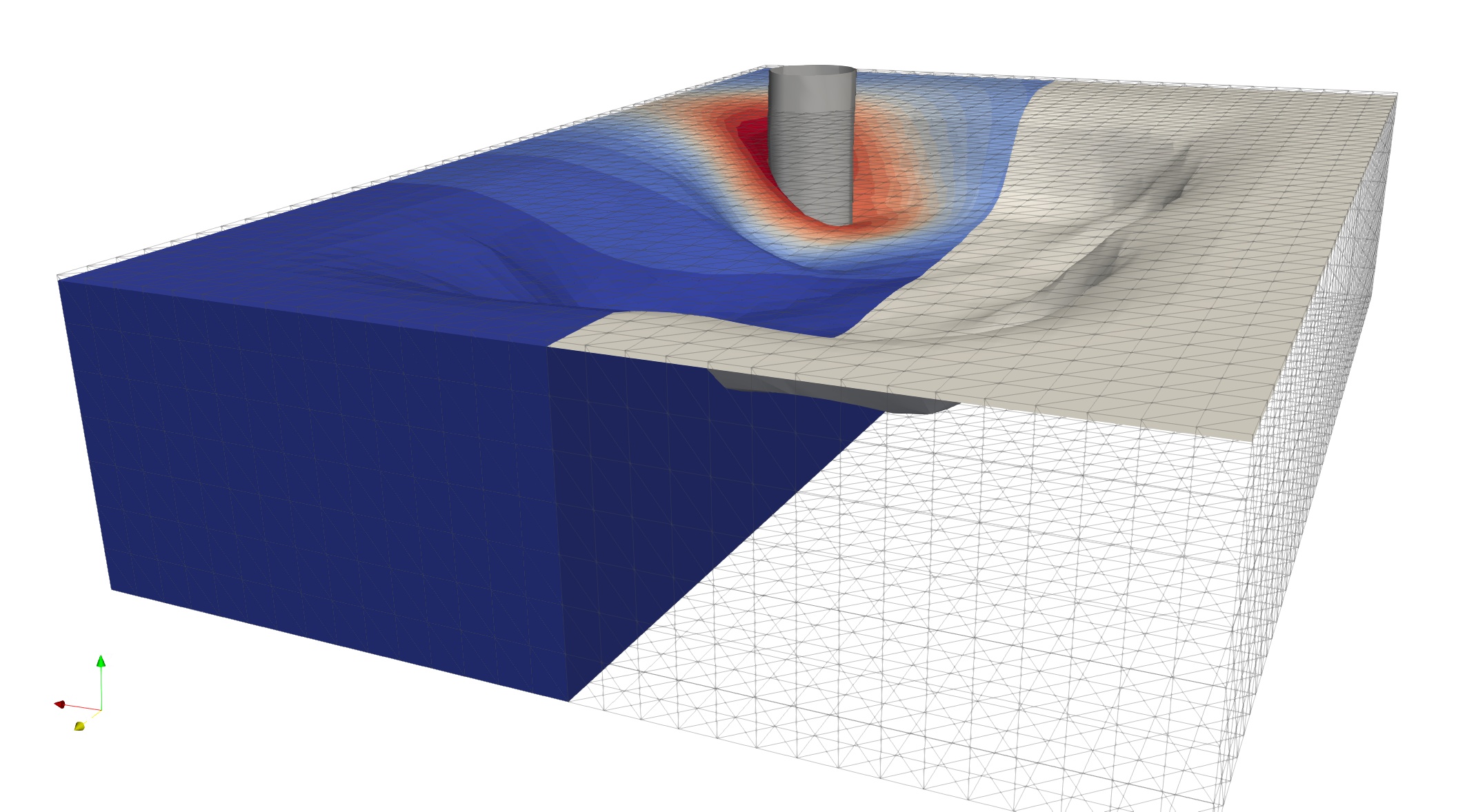}}\subfloat{\includegraphics[width=.1\textwidth]{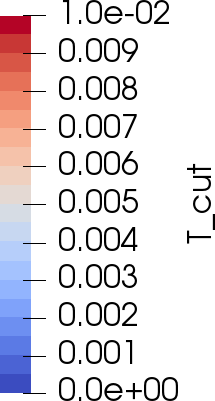}}
\end{center}
\caption{Laser beam going through the machining path shown in Figure~\ref{fig: 3D path} to create the desired rectangular cavity. The results correspond to analysis times $t=[0,1,2,3]$.}
\label{fig: temperature laser path 3D}
\end{figure}

\paragraph{Complex 3D machining path} We consider a complex machining path specified as shown in Figure~\ref{fig: 3D path}. We aim to form a rectangular cavity, using a machining strategy that is typical of what could be generated by a CAM software featuring thermal milling capabilities. The background domain is set to $\Omega_b = (-1,1) \times (-1.5,1.5) \times (-0.6,0.01)$ with an initial level set of $\phi(x,0) = z$. We set the time step to  $\Delta t =0.005$ and the mesh size to $h = 0.06$. The ablation strategy is described through the motion of the focal point $F(t)$ in time interval $t \in [0,3]$. The top layer is machined first, and deeper layers as represented in Figure~\ref{fig: 3D path} are applied subsequently. The remaining parameters of the laser beam are set to   $e_{ray} = (0,0,-1)$, $A_{amp}=3$ and $\sigma= 0.1$. We choose the material parameters as $T_m = 0.01$, $L=1$, $\rho=1$, $k=1$, $c=1$. For this particular example, function $f_t(t)$ is always equal to one (i.e. the laser fires continuously).
As shown in Figure~\ref{fig: temperature laser path 3D}, the manufacturing process  creates the expected rectangular cavity. The cylinder displayed in Figure~\ref{fig: temperature laser path 3D} represents the contour line $I=2$ of the laser beam. As the workpiece receives energy in a continuous way, no crater is formed. However, we can clearly see the streaks left by the laser beam, owing to a rather large hatch distance (i.e. distance between two consecutive straight lines of the machining path, within one particular layer). 

The results of this simulation can be played by paraview, using the .vtk files archived on zenodo.org by the authors \cite{claus_susanne_2018_1213279}. 



\section{Conclusions}

We have presented the first CutFEM algorithm dedicated to the solution of unsteady, one-phase Stefan-Signorini problems. The geometry of the domain is represented implicitly through the negative values of a continuous, piecewise linear level-set function defined using a regular, fixed finite element mesh. The boundary of the thermally ablated material can move arbitrarily and cut through the bulk of the elements, which circumvents the need for any remeshing operation during the simulation of phase change. We showed that the primal/dual formulation of the one-phase thermal ablation problem could be reformulated as a purely primal, nonlinear problem, using the Nitsche-Signorini idea, which avoids the need to introduce a Lagrange multiplier field for the interface velocity, and circumvents the need to design an inf-sup stable primal/dual discretisation strategy. Through the addition of stabilisation terms associated with the cut region, we proved that the method remains stable independently of the cut location. In addition, by carefully h-weighting several terms of the weak form associated with the proposed Stefan-Signorini-Nitsche method, we obtained optimal convergence with respect to spatial and temporal refinement. As a further contribution, we developed a 2D benchmark to test the convergence of numerical methods for one-phase Stefan problems. We hope that the new manufactured solution will be of use to researchers in the future. 

The robustness and versatility of the proposed algorithm was demonstrated through several representative examples in 2D and 3D. Although the method is general, our example section targeted realistic applications in laser micro-manufacturing, including the simulation of laser drilling and laser milling operations.


\clearpage

\section*{Acknowledgements}
The authors gratefully acknowledge the financial support provided by the Welsh Government and Higher Education Funding Council for Wales through the S\^{e}r Cymru National Research Network in Advanced Engineering and Materials under grants NRNG06 and NRN102.
\bibliographystyle{plainnat}
\bibliography{bibliography}

\end{document}